\documentclass[review,hidelinks,onefignum,onetabnum]{siamart250211}

\usepackage{amssymb}
\usepackage{appendix}
\usepackage{comment}
\usepackage{amsmath}
\usepackage{physics} 
\usepackage{mathtools}
\usepackage{bm}
\usepackage{booktabs}   
\usepackage{tabularx}   
\usepackage{makecell}   
\usepackage{xcolor}     
\usepackage{pifont}     
\usepackage{subcaption}

\usepackage{graphicx}    
\usepackage{multirow}  
\usepackage{geometry}
\usepackage[textsize=tiny]{todonotes}
\usepackage{xcolor}

\DeclareMathOperator*{\argmin}{arg\,min}

\newcommand{\cmark}{{\color{green!60!black}\ding{51}}} 
\newcommand{\xmark}{{\color{red!70!black}\ding{55}}}   

\renewcommand{\arraystretch}{1.25}

\usepackage{lipsum}
\usepackage{amsfonts}
\usepackage{graphicx}
\usepackage{epstopdf}
\usepackage{algorithmic}
\ifpdf
  \DeclareGraphicsExtensions{.eps,.pdf,.png,.jpg}
\else
  \DeclareGraphicsExtensions{.eps}
\fi

\newsiamremark{remark}{Remark}
\newsiamremark{hypothesis}{Hypothesis}
\crefname{hypothesis}{Hypothesis}{Hypotheses}
\newsiamthm{claim}{Claim}
\newsiamremark{fact}{Fact}
\crefname{fact}{Fact}{Facts}

\headers{Bilevel learning}{R. Grazzi, M. Pontil, S. Salzo, and A. Zemkoho}

\title{Bilevel learning\thanks{Submitted to the editors DATE.
\funding{This work was funded by EPSRC through a grant with reference EP/X040909/1, and in part by the European Union, under Projects Numbers 101070617 and 101120237.
}}}
\author{\scriptsize{Riccardo Grazzi}\thanks{Microsoft Research Cambridge, Cambridge, UK
  (\email{riccardo.grazzi.4@gmail.com}).}
\and \scriptsize{Massimiliano Pontil}\thanks{Computational Statistics and Machine Learning, Italian Institute of Technology, Genova, Italy \& AI Centre, Department of Computer Science, UCL, London, UK  
  (\email{m.pontil@cs.ucl.ac.uk}).}
\and \scriptsize{Saverio Salzo}\thanks{DIAG, Sapienza University of Rome, Via Ariosto, 25, 00185 Roma, Italy \&
Computational Statistics and Machine Learning, Italian Institute of Technology, Genova, Italy  (\email{salzo@diag.uniroma1.it}).}
\and \scriptsize{Alain Zemkoho}\thanks{School of Mathematical Sciences, University of Southampton, UK 
  (\email{a.b.zemkoho@soton.ac.uk}).}}

\usepackage{amsopn}

\ifpdf
\hypersetup{
  pdftitle={Bilevel Learning},
  pdfauthor={Riccardo Grazzi, Massimiliano Pontil, Saverio Salzo, and Alain Zemkoho}
}
\fi

\usepackage{amsmath, amssymb}
\usepackage{booktabs}

\begin{document}

\maketitle
\tableofcontents

\begin{abstract} 
Bilevel learning refers to machine learning problems that can be formulated as bilevel optimization models, where decisions are organized in a hierarchical structure. This paradigm has recently gained considerable attention in machine learning, as gradient-based algorithms built on the implicit function reformulation have enabled the computation of large-scale problems involving possibly millions of variables.
Despite these advances, the implicit function framework relies on restrictive assumptions, notably the requirement that the lower-level problem admit a unique optimal solution for each upper-level decision. Moreover, the computation of the derivative of the lower-level optimal solution function becomes significantly more involved when the lower-level problem includes constraints. As a result, many existing bilevel learning algorithms are effective only for relatively narrow classes of problems.
This paper reviews the main algorithmic ideas underlying recent progress in bilevel learning, highlighting both the key mechanisms responsible for their scalability and the limitations that arise in more general settings. We then draw connections with the broader bilevel optimization literature and discuss algorithmic techniques that may help overcome these limitations. Our aim is to bridge the gap between bilevel learning and classical bilevel optimization, thereby supporting the development of scalable methods capable of solving more general large-scale bilevel programs. 
\end{abstract}



\section{Introduction}\label{sec:Intro}

\textit{Bilevel learning} (BL) is a research area that lies at the intersection of \textit{bilevel optimization} (BO) and \textit{machine learning}. It focuses on machine learning problems with an inherent hierarchical structure that can be naturally modeled using the \textit{Stackelberg game} framework. In this setting, two decision makers interact: a leader (or \textit{upper-level player}) and a follower (or \textit{lower-level player}). The defining feature of a Stackelberg game is its sequential decision structure: the leader acts first, and the follower subsequently responds after observing the leader’s decision. This sequential order of play distinguishes Stackelberg games from other game-theoretic models, such as Nash or zero-sum games, in which decisions are made simultaneously. Consequently, Stackelberg games exhibit a vertical hierarchy between the players, with the leader occupying the upper level of the decision process. For this reason, they are often referred to as \textit{hierarchical games}, in contrast to the horizontal interaction that characterizes Nash-type games. A detailed introduction to the Stackelberg framework can be found in Heinrich von Stackelberg’s habilitation thesis \cite{Stackelberg1934} (see also the English translation \cite{von2010market}), as well as in the monograph \cite{dempe2002foundations}.

To avoid early distraction by the \textit{jungle} of bilevel optimization models, we defer a formal mathematical definition of a bilevel program to Section \ref{sec:What is a bilevel}, where different formulations of the problem and related concepts are presented. For the moment, we focus on the main purpose of this paper. When necessary, we distinguish between BO and BL. The motivation for this distinction is that BL is increasingly emerging as a research field in its own right and can therefore be characterized by the tools, approaches, and questions that arise in the design and analysis of solution algorithms. In particular, while BL is largely defined by the \textit{scale} of the problems it addresses, the traditional BO literature has been primarily concerned with the mathematical correctness of methods and theoretical properties of solution concepts. As a result, many BO approaches are tractable only for relatively small problem instances.

A key distinction between BL and BO therefore lies in the scale of the problems considered. Although many state-of-the-art BL methods are rooted in classical BO frameworks, the problems arising in BL can involve extremely large numbers of variables. For instance, the neural architecture design problem is modeled and solved in \cite{liu2019darts} as a bilevel optimization problem, where the resulting algorithm is applied to the CIFAR-10 and ImageNet datasets involving millions of parameters. Similarly, \cite{lorraine2020optimizing} develops a gradient-based algorithm for bilevel hyperparameter optimization in machine learning training and demonstrates its effectiveness on problems involving millions of weights and hyperparameters.

In contrast, many algorithms in the current BO literature remain largely conceptual or are designed primarily for small-scale problems (see, e.g., \cite{dempe2020bilevel}). Consequently, their direct applicability to large-scale BL settings is often limited. The \textit{first objective} of this paper is therefore to highlight, for the BO community, the algorithmic machinery that has recently enabled implicit-function–based gradient methods to scale to machine learning problems involving very large numbers of variables, as illustrated in the examples above. Note that the implicit function model, which underlies many state-of-the-art BL methods, consists essentially of substituting the lower-level optimal solution function, when it is well-defined as a vector-valued function, into the upper-level objective function.  This single-level reformulation then enables the development of gradient-based algorithms. 

It is worth recalling that the implicit function approach was already introduced in the early 1990s as a method for solving BO problems. In fact, it forms the basis of the principal solution approaches presented in two classical monographs on bilevel optimization, namely \cite{dempe2002foundations,outrata1998nonsmooth}. However, this framework did not gain widespread practical adoption, largely due to the difficulties involved in guaranteeing the existence of the lower-level optimal solution function and computing its Jacobian when it exists; see Sections \ref{sec:Main algorithmic techniques} and \ref{sec:Challenges and limitations} for further discussion. What has recently enabled the implicit function model to become a mainstream approach in BL is the introduction of derivative approximation approaches, inspired by advances in \textit{automatic differentiation} techniques, tools that have played a central role in the development of modern deep learning methods; see, e.g., \cite{baydin2018automatic,griewank2008evaluating}. 

Beyond derivative approximation techniques, another major departure from classical BO methods is the use of approximate solutions of the lower-level problem. In traditional BO approaches, the lower-level problem is typically assumed to be solved exactly, often to global optimality. In contrast, BL methods frequently rely on approximate solutions obtained through iterative optimization procedures. The recent resurgence of the implicit function framework in BL therefore comes with certain modeling simplifications. In particular, lower-level constraints are often avoided in order to facilitate the computation of gradient descent directions. However, many bilevel learning problems cannot realistically be formulated with unconstrained lower-level problems; see, e.g., \cite{ye2023difference,gao2023moreau,liu2023value,yao2024constrained,xu2023efficient}.

Moreover, the implicit function reformulation itself can be conceptually problematic, since the lower-level problem can typically admit multiple optimal solutions for some upper-level variable(s). Taking this issue into account—together with the fact that computing the Jacobian of the lower-level optimal solution function typically requires second-order information—the \textit{lower-level value function reformulation} has recently emerged as an important alternative in the BL literature (see, e.g., \cite{liu2022bome,kwon2023fully,kwon2023penalty,gao2022value,ye2023difference}). A key advantage of this approach is that it requires only first-order information to compute descent directions within gradient-based schemes, as it will be outlined in Section \ref{sec:Pure first and second order methods}.

The \textit{second objective} of this paper is to contribute to the acceleration of such developments by providing researchers in BL with a broad overview of algorithmic techniques from the BO literature. By bridging these two bodies of work, we aim to facilitate the development of enhanced tools capable of efficiently addressing large-scale BL problems that remain beyond the reach of existing approaches. More broadly, the goal of this article is to propose a unified perspective that brings together key algorithmic ideas from both BL and BO, thereby fostering the design of efficient algorithms for very large-scale problems and promoting new theoretical insights into the numerical behavior of bilevel programs.

\subsection{Related work}
A number of survey and overview papers on bilevel optimization and related topics already exist in the literature. On the BO side, to the best of our knowledge, the most recent and comprehensive collection of surveys is the edited volume \cite{dempe2020bilevel}. The topics covered in that volume are intentionally broad and include connections with other problem classes such as game theory and multiobjective optimization, as well as general theoretical questions and algorithms designed for particular classes of bilevel optimization problems. In contrast, the BO perspective adopted in this paper focuses primarily on algorithms for continuous bilevel optimization problems, while exploring the opportunities and challenges that arise when attempting to adapt or scale these methods to bilevel learning settings.
In addition, several earlier BO surveys and bibliographic reviews have provided high-level overviews of algorithmic developments for particular classes of bilevel optimization problems; see, for example, \cite{sinha2017review,Colson2007AnOO,kleinert2021survey,dempe2003annotated,vicente1994bilevel}.
On the BL side, a number of survey papers have appeared more recently. These works typically focus on specific machine learning applications of bilevel optimization and on tailored variants of gradient-based algorithms, together with empirical studies of their performance characteristics; see, for instance, \cite{liu2021investigating,cina2023wild,crockett2022bilevel,chen2022gradient,franceschi2025}. A more recent overview aimed at a broader machine learning audience is provided in \cite{zhang2024introduction}, which explicitly discusses the connections between bilevel optimization and modern signal processing and machine learning applications.

\subsection{Main contributions}
In summary, the main contributions of this paper are as follows:
(i) to provide an overview of state-of-the-art numerical schemes tailored to BL that have been developed to address medium- to large-scale problems, while highlighting the key features that drive their efficiency;
(ii) to identify the limitations and weaknesses of existing BL algorithms, thereby enabling researchers in BO interested in machine learning applications to better understand the challenges and shortcomings of current approaches;
(iii) to present a collection of state-of-the-art algorithmic techniques from the broader BO literature for general bilevel programs, allowing BL researchers to quickly familiarize themselves with potential methodological directions for addressing problems that current BL algorithms cannot effectively solve; and
(iv) to provide a reference resource for researchers interested in (or beginning work on) bilevel optimization and its applications in machine learning.

\subsection{Organization of the paper}
In the next section, we introduce the mathematical description of the \textit{bilevel optimization} concept and its different formulations. This step is necessary because the term bilevel optimization is used to describe a variety of problems that may differ in structure while still preserving the defining characteristics of a Stackelberg game model, as outlined above. It is therefore useful to begin by presenting the main variants of the problem and clarifying how they are related and how they build upon one another.
Subsequently, in Section \ref{sec:A short history}, we provide a brief historical overview of BO and illustrate how many problems arising in BL naturally exhibit a Stackelberg structure, even though this connection was not always explicitly recognized or formally linked to the BO literature.

For readers who are less familiar with applications of BO in machine learning, Section \ref{sec:A flavor of machine learning} presents a concise overview of bilevel programming applications in machine learning (i.e., an overview of \textit{BL problems}). We also discuss in greater detail two representative examples—hyperparameter optimization and adversarial learning—which correspond to special classes of optimistic and pessimistic bilevel programs, respectively (these concepts will be formally defined in Section \ref{sec:What is a bilevel}).
Section \ref{sec:Main algorithmic techniques} focuses on the state-of-the-art methods currently used in BL, which are largely built upon the implicit function reformulation. This framework relies on the well-posedness of the lower-level problem, typically requiring the existence of a unique optimal solution for every upper-level variable. We discuss the main assumptions underlying this approach—such as strong convexity and the absence of lower-level constraints—as well as the algorithmic mechanisms that have enabled its practical success, including approximation strategies for derivatives and lower-level solutions. Convergence paradigms and performance evaluation techniques are also reviewed in this section.

In Section \ref{sec:Challenges and limitations}, we examine several important limitations and challenges associated with the state-of-the-art methods described in Section \ref{sec:Main algorithmic techniques}, and discuss possible directions for overcoming these difficulties.
Section \ref{sec:Constrained optimization} then provides an overview of algorithmic techniques from the BO literature that can be used to address problems in which the lower-level problem is convex—but not necessarily strongly convex—and may include constraints. Subsequently, in Section \ref{sec:Pure first and second order methods}, we relax the convexity assumption on the lower-level problem and explore potential avenues for developing algorithms for BO that rely purely on first- or second-order information. Such approaches are not directly applicable within the framework discussed in Section \ref{sec:Main algorithmic techniques}, since the Jacobian of the lower-level optimal solution function typically requires second-order information about the lower-level problem.

Overall, the study in this paper is mainly focused on the \textit{standard optimistic} bilevel optimization problem (see next section for the definition), and the algorithmic discussions throughout Sections \ref{sec:Main algorithmic techniques}--\ref{sec:Pure first and second order methods} are dedicated to different transformations of the model; namely, the implicit function, Karush-Kuhn-Tucker, and lower-level value reformulations. In Section \ref{sec:Comparing the reformulations}, a comparison of these different approaches is provided. 
Finally, Section \ref{sec:Some final thoughts} concludes the paper with a set of observations and perspectives on how the ideas discussed in Sections \ref{sec:Constrained optimization} and \ref{sec:Pure first and second order methods} could be further explored in the context of BL, where they remain largely underdeveloped.

\section{What is a bilevel optimization problem?}\label{sec:What is a bilevel}

This is possibly one of the most fascinating questions, as there are multiple concepts labeled as/or related to \textit{bilevel optimization}, depending on what interpretation one makes, or also the specific area of application. The first and most widely understood, and which will also be at the center of the attention of this paper, because it is the one generally used in machine learning, is the problem 
\begin{equation}\label{eq:Bilevel-Optimization-Problem}\tag{BOP}
      \underset{x\in X}{\min}~F(x, y)\;\;  \mbox{s.t.}\;\;   y\in S(x):=\underset{y\in Y(x)}{\text{argmin}}~f(x,y),
\end{equation}
where the functions $F : \mathbb{R}^n\times \mathbb{R}^m \rightarrow \mathbb{R}$ and $f : \mathbb{R}^n\times \mathbb{R}^m \rightarrow \mathbb{R}$ represent the upper- and lower-level objective functions, respectively. Similarly, $X \subset  \mathbb{R}^n$ is the upper-level feasible set, while the set-valued mapping $Y: X \rightrightarrows \mathbb{R}^m$ describes the lower-level feasible set. Overall, \eqref{eq:Bilevel-Optimization-Problem} represents the upper-level (or leader's) problem, while the set $S(x)$ collects all the optimal solutions of the lower-level (or follower's) problem (for all $x\in X$): 
\begin{equation}\label{eq:Lower-level-problem}\tag{LL}
     \min_y~f(x,y) \;\;\mbox{ s.t. }\;\; y\in Y(x).
\end{equation}

Problem \eqref{eq:Bilevel-Optimization-Problem} will be said to be well-posed if one assumes that for any choice $x\in X$ of the leader, the follower has a single optimal  solution. Precisely, this means that the following condition, where $|C|$ stands for the cardinality of the set $C$, is satisfied:
\begin{equation}\label{eq:S(x)=1}
    \left\{x\in X:\;\; |S(x)|= 1 \right\} = X;
\end{equation}
i.e., we have $S(x)=\{y(x)\}$ for all $x\in X$ with $y(\cdot): X \rightarrow \mathbb{R}^m$ being the optimal solution function of the lower-level problem. In this case, problem \eqref{eq:Bilevel-Optimization-Problem} reduces to
\begin{equation}\label{eq:Implicit-Function-Model}
\tag{\mbox{$\text{P}_{\text{i}}$}}
      \underset{x\in X}{\min}\; \mathcal{F}(x) := F(x, y(x)).
\end{equation}
This model, known in the classical bilevel optimization literature as \textit{implicit function reformulation} (see \cite{kolstad1990derivative} for one of the very first studies based on the approach) is the state of the art working framework in  \textit{bilevel learning}.

However, there is also another very rich set of solution concepts for  \eqref{eq:Bilevel-Optimization-Problem} when the lower-level optimal solution set-valued mapping $S$ satisfies the condition
\begin{equation}\label{eq:S(x)>1}
    \left\{x\in X:\;\; |S(x)| >1 \right\} \neq \emptyset.
\end{equation}
This corresponds to the situation where the lower-level problem has more than one optimal solution for some choices of the upper-level player. In the context of \eqref{eq:S(x)>1}, there are two radically opposed solution concepts for problem \eqref{eq:Bilevel-Optimization-Problem}. 
The first, and more commonly used one, is to assume that each time the leader picks an $x\in X$ that leads to $|S(x)|>1$, the follower selects a value $y\in S(x)$ that is in favorable to the leader. This leads to  the (original) \textit{optimistic} bilevel optimization problem 
\begin{equation}\label{eq:Optimistic-Model}\tag{\mbox{$\text{P}_{\text{o}}$}}
      \underset{x\in X}{\min}~\underset{y\in S(x)}{\min}~F(x, y),
\end{equation}
also known as the cooperative model. The optimistic bilevel program is the most studied class of the problem; however, investigations are almost never on the formulation \eqref{eq:Optimistic-Model}, but rather on the following version of the problem, labeled in \cite{DempeMordukhovichZemkoho2012sensitivity,zemkoho2016solving} as  \textit{standard} optimistic bilevel optimization problem in opposition to the original optimistic \eqref{eq:Optimistic-Model}: 
\begin{equation}\label{eq:Standard-Optimistic-Bilevel}\tag{P}
      \underset{x,\, y}{\min}~F(x, y)\;\;  \mbox{s.t.}\;\; x\in X, \;\;  y\in S(x).
\end{equation}
It can be seen that here, full control over the leader and follower's variables $x$ and $y$ is given to the upper-level player. 
Although it might sound intuitive that problems \eqref{eq:Optimistic-Model} and 
\eqref{eq:Standard-Optimistic-Bilevel} are equivalent, it was shown in \cite{DempeMordukhovichZemkoho2012sensitivity} that this is true only if global optimal solutions are considered. But locally, both problems are not equivalent. For a local optimal solution $\bar x$ of problem \eqref{eq:Optimistic-Model}, any point $(\bar x, \bar y)$ with $\bar y\in S(\bar x)$ is a local optimal solution of problem \eqref{eq:Standard-Optimistic-Bilevel}. However, if $(\bar x, \bar y)$ is locally optimal for problem \eqref{eq:Standard-Optimistic-Bilevel}, one needs the set-valued mapping
\[
S_{\text{o}}(x):= \underset{\quad\;\, y\in S(x)}{\text{argmin}}~F(x,y)
\]
to be \textit{inner semicontinuous} at $(\bar x, \bar y)$ to ensure that $\bar x$ is locally optimal for problem \eqref{eq:Optimistic-Model}. This is quite a strong assumption; for its definition, and more details on the relationship between these two problems, interested readers are referred to \cite{DempeMordukhovichZemkoho2012sensitivity}.

The second option, which is a bit less researched, is to assume that each time the leader picks an $x\in X$ such that  $S(x)>1$, the follower selects a value $y\in S(x)$ that is antagonistic to the leader. To anticipate on unfavorable choices from the follower, the leader solves the so-called \textit{pessimistic} bilevel optimization problem 
\begin{equation}\label{eq:Pessimistic-Model}\tag{\mbox{$\text{P}_{\text{p}}$}}
      \underset{x\in X}{\min}~\underset{y\in S(x)}{\max}~F(x, y).
\end{equation}
This is a more challenging problem to solve, and in the hope of promoting the implementation of ideas from the abundant literature on the standard optimistic bilevel program  \eqref{eq:Standard-Optimistic-Bilevel}, the  \textit{standard pessimistic} bilevel optimization problem
\begin{equation}\label{eq:Standard-Pessimistic-Bilevel}
      \underset{x,\, y}{\min}~F(x, y)\;\;  \mbox{s.t.}\;\; x\in X, \;\;  y\in S_{\text{p}}(x)
\end{equation}
was recently investigated in \cite{lampariello2019standard}. Note that here, the \textit{two-level} optimal solution set-valued mapping $S_{\text{p}}: X \rightrightarrows \mathbb{R}^m$ is defined by
\[
S_{\text{p}}(x):= \underset{\quad\;\, y\in S(x)}{\text{argmax}}~F(x,y).
\]
A combination of  \textit{Karush-Kuhn-Tucker} and some \textit{lower-level value function}--type reformulation (see Section \ref{sec:Constrained optimization} and Section \ref{sec:Pure first and second order methods} for some relevant details) is used to address $S_{\text{p}}$ in order to get a single--level transformation of problem \eqref{eq:Standard-Pessimistic-Bilevel}.

The optimistic and pessimistic bilevel optimization problems \eqref{eq:Optimistic-Model} and \eqref{eq:Pessimistic-Model}, respectively, represent two very extreme positions, either the leader and follower cooperate (optimistic problem) or they do not (pessimistic problem). For this reason, a number of works (see, e.g., \cite{cao2002partial,lagos2023complexity}) have suggested a compromise model, where in a nutshell, if the follower has multiple options to pick from, for some choices of $x\in X$, they select one that represents a compromise between the leader and the follower; i.e., the problem to be solved is 
\begin{equation}\label{eq:Partial-Cooperation-Model}\tag{\mbox{$\text{P}_{\text{op}}$}}
      \underset{x\in X}{\min}~\lambda\varphi_{\text{o}}(x) + (1-\lambda)\varphi_{\text{p}}(x),
\end{equation}
with $\lambda \in [0, \, 1]$ denoting the degree of cooperation between the upper- and lower-level players \cite{shihui2013new}. In problem \eqref{eq:Partial-Cooperation-Model}, the functions $\varphi_{\text{o}}$ and $\varphi_{\text{p}}$ are respectively defined by 
\begin{equation}
  \varphi_{\text{o}}(x) := \min\,\left\{F(x, y):\, y\in S(x)\right\} \mbox{ and } \varphi_{\text{p}}(x) := \max\,\left\{F(x, y):\, y\in S(x)\right\}. 
\end{equation}
Labeled as \textit{two-level} optimal value functions and studied in detail in \cite{DempeMordukhovichZemkoho2012sensitivity}, where suitable conditions for their local Lipschitz conditions are established. 

Finally, \cite{aboussoror1995strong,aboussoror2017strong} introduces the following \textit{weak-strong} \textit{Stackelberg}/\textit{bilevel optimization} problem as a generalization of the above optimistic and pessimistic bilevel programs: 
\begin{equation}\label{eq:weak-strong-model}
   \underset{x\in X}{\min}~\underset{y\in S(x)}{\min}~\underset{z\in S(x)}{\max}~\mathfrak{F}(x, y, z). 
\end{equation}
Here, the function $\mathfrak{F} : \mathbb{R}^n\times \mathbb{R}^m \times \mathbb{R}^m \rightarrow \mathbb{R}$ represents the upper-level objective function, while $S: X \rightrightarrows \mathbb{R}^m$ describes the lower-level optimal solution set-valued mapping defined in \eqref{eq:Bilevel-Optimization-Problem}. It might be useful to recall that the optimistic (resp. pessimistic) bilevel program \eqref{eq:Optimistic-Model} (resp. \eqref{eq:Pessimistic-Model}) is also called \textit{weak} (resp. \textit{strong}) Stackelberg/bilevel optimization problem. Hence, the reason why the problem is referred to as weak-strong Stackelberg problem. Similarly, it could therefore be labeled as \textit{optimistic-pessimistic} bilevel optimization problem. To see why this vocabulary makes sense, observe that if $\mathfrak{F}(x, y, z) := F(x, y)$, we get the original optimistic model \eqref{eq:Optimistic-Model}, while having $\mathfrak{F}(x, y, z) := F(x, z)$ leads to the pessimistic bilevel optimization problem \eqref{eq:Pessimistic-Model}. Moreover, setting $\mathfrak{F}(x, y, z) := \lambda F(x, y) + (1-\lambda) F(x, z)$, we can observe that we get the partial cooperation model in \eqref{eq:Partial-Cooperation-Model} for a fixed $\lambda\in [0,\, 1]$.

\section{A short (possibly) shared history}\label{sec:A short history}
Bilevel optimization emerged from the habilitation thesis of von Stackelberg in 1934 \cite{Stackelberg1934}. First, it attracted interest mainly from economists (see, e.g.,  \cite{LeonidHurwicz1945,Kornai1965,ChenCruz1972}; for a thorough economic perspective on von Stackelberg's work, see {\href{https://www.emerald.com/insight/publication/issn/0144-3585/vol/23/iss/5/6}{Volume 23 Issue 5/6 of the Journal of Economic Studies} specifically dedicated to his research contribution and its influence) until 1973 when it was introduced to the field of mathematical optimization by Bracken and McGill \cite{bracken1973mathematical}. It is however important to note that the problem introduced in Bracken and McGill's first paper on the subject instead corresponds to what is known today as a semi-infinite programming problem \cite{hettich1993semi}.
Candler and Norton in their  reports \cite{candler1977multiPolicy} and \cite{candler1977multi} might have been the first to connect the dots between the bilevel optimization model that emerged in the mathematical optimization literature and the work of von Stackelberg, and can also be credited for coining the expression ``bilevel optimization'' (namely, multilevel for optimization problems with more than two levels). 

The pioneering works of operations research and mathematical optimization experts from around the early 70s came with two important things: (i)  the emergence of applications of bilevel optimization outside of economics and pure game theory; for example, the initial works of Bracken and McGill focused on military and defense applications \cite{bracken1973mathematical,bracken1974defense}, while Candler and Norton highlighted the applicability of the problem in areas such as engineering, biology, and policy design and implementation in agriculture  \cite{candler1977multiPolicy,candler1977multi}. 
(ii) A huge interest in the development of solution algorithms to solve bilevel optimization problems. 
Works on solution algorithms started to intensify around the early 1980s, as it can be seen in this first survey on the subject by Kolstad \cite{kolstad1985review} in 1985. 
Overall, the algorithms developed so far could be generally classified in the four categories outlined in following subsections. 

\subsection{Methods for bilevel programs with fully linear functions} The problem \eqref{eq:Bilevel-Optimization-Problem} will be said to be fully linear if the functions $F$ and $f$ are linear in $(x,y)$ and the sets $X$ and $Y(x)$, for all $x\in X$, are defined by functions that are linear in $x$ and $(x,y)$, respectively. Problem \eqref{eq:Standard-Optimistic-Bilevel} in this case has an interesting geometric structure, in the sense that at least one optimal solution of the problem, if one exists, occurs at a vertex of the polyhedron the set $\left(X\times \mathbb{R}^m\right) \cap \text{gph}Y$. Note that here and in the sequel, for a set-valued mapping $\Gamma: \mathbb{R}^n \rightrightarrows \mathbb{R}^m$, its \textit{graph} is defined by 
    \[
    \text{gph}\,\Gamma:=\left\{(x, y)\in \mathbb{R}^n\times \mathbb{R}^m\left|~y\in \Gamma(x)\right.\right\}.
    \]
    This result, first discovered by Bials and Karwan \cite{BialasKarwan1982}, has been the basis for the development of multiple enumerative algorithms for fully linear bilevel programs; see, e.g., \cite{CANDLER198259,Papavassilopoulos1982,NarulaNwosu1983}. Unfortunately, as we will see in the general overview of bilevel learning problems, it is not clear that there are any applications of a fully linear model in the context machine learning, except perhaps that it might happen that for some problems, such cases could appear as subproblems for some algorithmic problems given that, for example, as it will be clear in Section \ref{sec:Main algorithmic techniques}, the process of computing the directional derivative for the lower-level problem can be seen as solving a quadratic linear bilevel program. 

    \subsection{Algorithms tailored to problems with convex lower-level problems}  The lower-level problem \eqref{eq:Lower-level-problem} is said to be convex if the objective function $f(x, .)$ and feasible set $Y(x)$ are both convex for all  $x\in X$. Of course, the fully linear case above is a special case of this one, but in this more general setting, the construction of numerical methods has generally relied on the so-called Karush-Kuhn-Tucker (KKT) reformulation, which replaces the lower-level problem in \eqref{eq:Standard-Optimistic-Bilevel} by its corresponding KKT conditions. If the lower-level problem has inequality constraints, this reformulation leads to a problem with complementarity constraints, which introduces a very high-level of complexity in the problems.
    
    Various methods have been proposed to solve this problem, with the main focus usually being on how to handle the complementary constraints; early related papers include the work of Fortuny-Amatand McCarl \cite{Fortuny-AmatMcCarl1981}, who introduce the famous big M method, which has been very influential in the field for many year, the work of Bard and Falk \cite{bard1982explicit}, which introduces a branch and bound--type method, which consists to solve convex approaximations of the KKT reformulation at each iteration, as well as the work by Bialas and Karwan \cite{bialas1982two}, where a pivot-type method is developed to compute an approximation of the KKT reformulation as a problem of finding the solution of a mixed-complementary system. As it will be shown in Section \ref{sec:Constrained optimization}, despite being fundamentally a nonconvex constrained optimization problem, the KKT reformulation has a potential for bilevel learning problems, as it is less restrictive in terms of the required framework. 

      \subsection{Techniques based on strongly convex lower-level problems} This is the framework that enables the implicit function model \eqref{eq:Implicit-Function-Model} to be well-defined, and has been the main based for the development of numerical algorithms in the context of BL. One of the main motivations of this paper has come from witnessing the huge interest that the implicit model \eqref{eq:Implicit-Function-Model} has attracted in the context of solving bilevel programs appearing in machine learning. The resurgence of methods based on the implicit has been a source of curiosity, as progress on the use of such methods seemed to have stalled in the more general field of BO. Considering the performance of this approach and the depth of analysis of the corresponding algorithms in BL, we aim to identifying the reasons of this success and draw attention to lessons that could be learned for other applications of BO. 
      Before we come back to the technicalities of such methods and reasons for their success in Section \ref{sec:Main algorithmic techniques}, we would like to point out that gradient descent approach, which has been the main algorithmic technique in this context, has been in existence in bilevel optimization since the early development of mathematical optimization--based numerical methods for the problem. 
      
      The PhD thesis of de Silva \cite{de1979sensitivity} completed under the supervision of Garth McCormick, who was at the forefront of the development of sensitivity analysis for optimal solutions of parametric optimization problem, was probably the first work in this context, implementing mainstream implicit function results for the calculation of $\nabla y(x)$ for a gradient descent scheme for a problem of the form \eqref{eq:Implicit-Function-Model}; see the paper \cite{desilva1992implicitly} with some of the related results. 
 Around the same period,  Shimizu and Aiyoshi \cite{shimizu1981new} propose another gradient descent scheme, where a barrier approach is used to eliminate lower-level constraints, before an implicit function technique is applied to the resulting Fermat rule of the new penalized problem. However, the work of Kolstad and Lasdon \cite{kolstad1990derivative}, well-known for a gradient descent scheme for problem \eqref{eq:Implicit-Function-Model}, is probably the first article that introduced an approximation approach for $\nabla y(x)$, which was then applied to solve relatively large size problems at the time (see details in  \cite{kolstad1986derivative}). 
    Bundle methods for \eqref{eq:Implicit-Function-Model} have also been very prominent in solving problem \eqref{eq:Implicit-Function-Model}, and represent  the main focus of the book \cite{outrata1998nonsmooth}; it was also  an important component of the exposition in other books on bilevel optimization such as  \cite{bard2013practical,dempe2002foundations}. 

       \subsection{General bilevel program without any explicit convexity requirement}
       In this case, the natural way to transform problem \eqref{eq:Standard-Optimistic-Bilevel} into a numerically tractable problem has been through the lower-level value function consisting to replace the inclusion $y\in S(x)$ by its definition. Initial ideas in this context can be traced back to \cite{bracken1974method,geoffrion1972coordination}. More recent effort in this include the works \cite{mitsos2008global,wiesemann2013pessimistic,kleniati2014branch,paulavivcius2016global}, which largely exploit the connection of the value function reformulation to semi-infinite programming to build methods to compute global optimal solutions. Considering the underlying techniques, the approaches in these papers are unlikely to scale well to the size of problems encountered in machine learning. This is because the schemes in the latter articles either rely on some branch-and-bound techniques (\cite{mitsos2008global,kleniati2014branch,paulavivcius2016global}) or semi-infinite programming--type discretization techniques (\cite{mitsos2008global,wiesemann2013pessimistic,kleniati2014branch,paulavivcius2016global}). A second class of method that has emerged recently in the literature (see \cite{zemkoho2021theoretical,fischer2022semismooth,fliege2021gauss}), building on standard continuous nonlinear optimization theory, and to be overviewed in Section \ref{sec:Pure first and second order methods}, has more potential in the context of bilevel learning. 


\subsection{Some bilevel learning history and connections to bilevel optimization} The BL history can go as far behind as well, if we consider the multiple problems that are now well-understood as bilevel optimization problems, but had stayed in the shadows of the field for a very long time. For example, in Bengio's paper \cite{bengio2000gradient}, where he studies the gradient--based approach for hyperparameter optimization in machine learning,  
which is simply an early version of the now classical gradient descent method for the implicit function model \eqref{eq:Implicit-Function-Model}, discussed in Section \ref{sec:Main algorithmic techniques} of the bielevel optimization problem, the Akaike Information Criterion (AIC) model is referred to as one of the long standing techniques for hyperparameter computation \cite{akaike1974new}. 
If we look closely at the AIC model in statistics, one of its applications is for cross-validation in the context of training the autoregressive integrated moving average (ARIMA) model for forecasting, it is used to find the order $(p, d, q)$, and this problem can be written as the bilevel program 
\begin{equation}\label{eq:ARIMA}
      \underset{p, d, q}{\min}~{\text{AIC}}(p, d, q, \theta, \phi)\;\;  \mbox{s.t.}\;\;  (\theta, \phi)\in \underset{\theta, \phi}{\text{argmin}}~f(p, d, q, \theta, \phi),
\end{equation}
where the order $(p, d, q)\in \mathbb{N}^3$ and $\theta$ and $\phi$ represent the AR and MA model parameters, respectively. The AIC model dates back to 1974 \cite{akaike1974new} and even if it has probably not yet been written in the form \eqref{eq:ARIMA}, it can naturally be translated as such. Similar BO models can be written to compute hyperparameters in a wide range of statistics problems, including design of experiments \cite{eriksson2000design}, regression analysis \cite{draper1998applied}, as well the estimation of parameters in various areas of engineering; see, e.g., \cite{mitsos2009bilevel,bollas2009bilevel,glass2018bilevel}. It might be worth mentioning here that hyperparameter optimization in machine learning is one of the most prominent areas of applications of BO in machine learning, as we will discuss in the next section. 

In terms of direct connections between machine learning and BO, it seems like links started to become clearer in the works of Bengio \cite{bengio2000gradient} and Chapelle et al. \cite{chapelle2002choosing}, where they developed gradient descent--type algorithms for hyperparameter computation based on the implicit function model (see Fig. \ref{fig:ChapelleEtal-Algorithm}), where the implicit function model is used, which is now the classical framework for bilevel learning algorithms, as it will be discussed in Section \ref{sec:Main algorithmic techniques}.

\begin{figure}
    \centering
    \includegraphics[scale = 0.45]{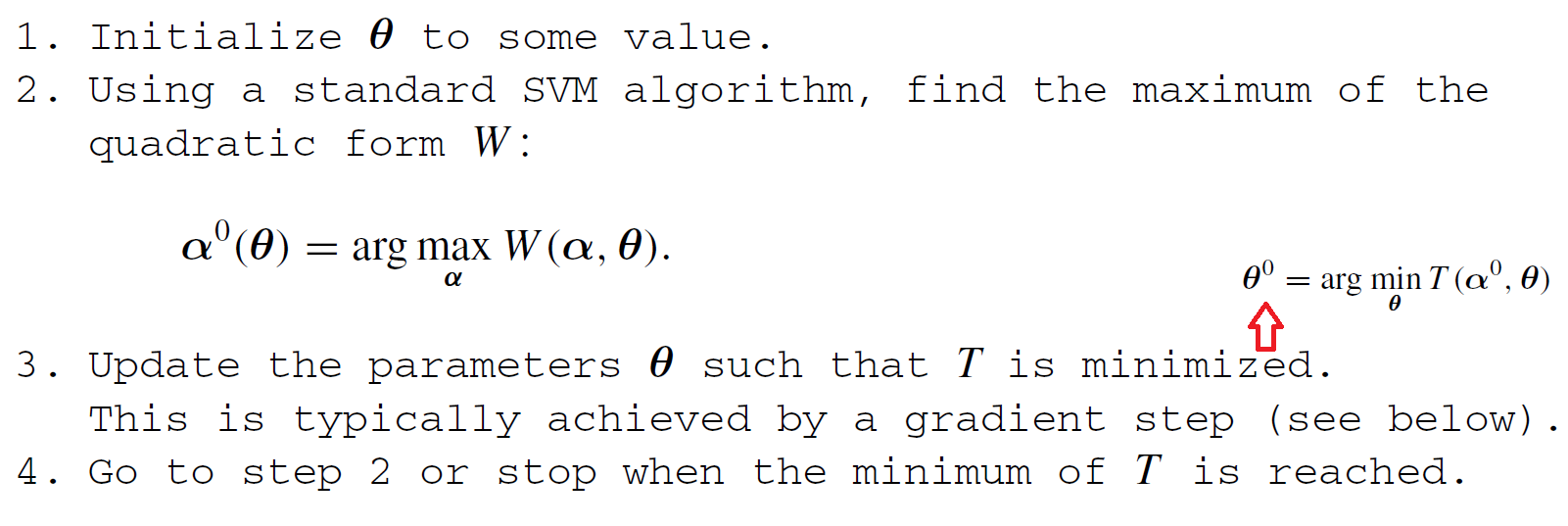}
    \caption{Algorithm for  bilevel hyperparameter optimization  by Chapelle et al. \cite{chapelle2002choosing}}
    \label{fig:ChapelleEtal-Algorithm}
\end{figure}

It must be said that in the articles \cite{bengio2000gradient} and Chapelle et al. \cite{chapelle2002choosing}, there is no mention of bilevel optimization or Stackelberg games; moreover, no relevant articles on the subject seems to be mentioned. So, their works could potentially be cast as independent discoveries of bilevel optimization. To the best of our knowledge, Kristin Bennett and her co-authors in a series of papers mainly on hyperparameter optimization for support vector machines \cite{bennett2006model,bennett2006interplay,Kunapuli2007BilevelMS,Bennett2008ABO,bennett2008bilevel,Kunapuli2008ClassificationMS,moore2009nonsmooth,Bennett2010BilevelPA} can be credited as the first to clearly make the connection between bilevel learning and bilevel optimization. Since the publication of this series of papers between 2006 and 2010, bilevel optimization has literally exploded in machine learning, and seems to be have become an area of research in its own right, within the field of machine learning. To illustrate this, we can mention, for example, the surveys \cite{liu2021investigating,crockett2022bilevel,chen2022gradient,zhang2023introduction}, just in the last 5 years, and a couple of hundreds of papers on the subject have been published in the last few years. So many applications of bilevel optimization have been discovered in machine learning; see next section for a sample of them.  

A key focus of this paper is to look closely at the technical aspects of the main BL algorithmic framework, in order to highlight what can be learned from it to tackle other bilevel programs. Conversely, as most algorithms in the current BL  literature are based on the implicit function model, this paper aims to draw the attention of BL researchers to the multitude of approaches that could be used to address a wide range of problems for which the model \eqref{eq:Implicit-Function-Model} is not well-defined, because of the failure of condition \eqref{eq:S(x)=1}.

\section{A flavor of bilevel optimization applications in machine learning}\label{sec:A flavor of machine learning}
There is a wide range of applications of BO in machine learning, and the number has been growing steadily in recent years; see, e.g., \cite{liu2021investigating,cina2023wild,crockett2022bilevel,chen2022gradient} for surveys focused around applications. Our aim here is to give a flavor of how typical such applications look like, while focusing on where the challenges of BL problems lie. Recall that the major difference that sets apart these applications from the ones usually considered in the classical BO literature is the scale of the considered problems. Indeed, the number of lower-level variables, the number of samples in the dataset, and in some cases even the number of upper-level variables can reach millions or even billions \cite{franceschi2018bilevel, lorraine2020optimizing, pan2025scalebio}. In this setting, scalability  becomes a core priority for the development of practically relevant algorithms.


We first focus our attention on two standard examples representing the extreme modeling viewpoints discussed in Section \ref{sec:What is a bilevel} occurring when condition \eqref{eq:S(x)>1} holds; i.e., the optimistic and pessimistic models, represented by the \textit{hyperparameter optimization} and \textit{data poisoning} problems, respectively. Let $\mathcal{D}_{\text{tr}}$ denote the \textit{training set} used to optimize the lower-level model parameters $y \in \mathbb{R}^m$ (e.g., neural network weights), and let $\mathcal{D}_{\text{val}}$ denote the \textit{validation set} used to evaluate the upper-level decision variables $x \in \mathbb{R}^n$ (e.g., hyperparameters).
In the supervised setting, the datasets consist of input-target pairs; for instance, $\mathcal{D}_{\text{val}} = \{\xi_i\}_{i=1}^{N_{val}}$ where each $\xi_i = (u_i, v_i)$ represents a data point $u_i$ and its label $v_i$ for $i=1, \ldots, N_{val}$.

\subsection{Hyperparameter optimization in machine learning}
In a hyperparameter optimization (HPO) problem, the goal is to select hyperparameters $x \in X$ that minimize a validation loss $F$ (upper-level objective function), subject to the model parameters $y$ minimizing a training loss $f$ (lower-level objective function). In this cooperative setting, the leader assumes the follower selects the $y \in S(x)$ most favorable to the upper-level objective. The optimistic formulation of the HPO problem reads as 
\begin{equation} \label{eq:Optimistic-Concrete}
    \min_{x \in X} \inf_{y \in S(x)} F(x, y) \coloneqq \sum_{\xi \in \mathcal{D}_{\text{val}}} \ell_{\text{val}}(y, x; \xi),
\end{equation}
where $S: X \rightrightarrows \mathbb{R}^m$ describes the optimal solution set-valued mapping of the training (or lower-level) problem; to be more precise, 
\begin{equation}\label{eq:lower-level-concrete}
     S(x) \coloneqq \argmin_{z \in \mathbb{R}^m}  f(x, z) \coloneqq \sum_{\zeta \in \mathcal{D}_{\text{tr}}} \ell_{\text{tr}}(z, x; \zeta) + \mathcal{R}(z, x).
\end{equation}
Here, $\mathcal{R}$ is a regularizer (e.g., $\ell_2$ regularization). 
Here, $\zeta \in \mathcal{D}_{\text{tr}}$ represents a training sample pair. The functions $\ell_{\text{tr}}$ and $\ell_{\text{val}}$ typically measure the discrepancy between the model prediction and the true label. Common choices for the validation loss $\ell_{\text{val}}$ include the Squared Error for regression tasks or the Cross-Entropy loss for classification. The key distinction is that $\ell_{\text{tr}}$ is evaluated on the training set and is often modulated by $x$. A prominent example is \textit{data reweighting} \cite{ren2018learning}, where $x$ assigns a weight to each training sample, leading to the expression $\ell_{\text{tr}}(z, x; \zeta_i) = \sigma(x_i) \ell(z; \zeta_i)$ (where $\sigma$ ensures positivity), allowing the leader to downweight noisy or corrupted data. 

The term $\mathcal{R}(z, x)$ represents a regularization term explicitly controlled by the hyperparameters. A standard example is weight decay, where $\mathcal{R}(z, x) = \sum_{i=1}^m \exp(x_i) z_i^2$, with $x$ acting as the log-regularization coefficient for each weight.
In practice, $X$ often includes continuous and discrete parameters such as \textit{regularization coefficients} and \textit{architecture depth}. Seminal works have applied this framework to kernel methods and Support Vector Machines (SVMs), where the lower-level problem often enjoys convexity \cite{chapelle2002choosing, keerthi2006efficient, bennett2008bilevel}. Other optimistic BL problems can be found here \cite{liu2022general,liu2020generic,sow2022primal,yao2024constrained,xiao2023generalized}. More general gradient-based approaches for continuous hyperparameters have been explored in \cite{pedregosa2016hyperparameter, franceschi2018bilevel}.  

Historically, hyperparameter optimization was predominantly performed using Grid Search or Random Search. While effective for a small number of hyperparameters (e.g., $n < 10$), these ``black-box'' methods scale exponentially poorly with dimension. The BO perspective represented a paradigm shift by utilizing gradient information from the lower-level (via implicit differentiation or iterative differentiation). This allows for the simultaneous optimization of thousands of hyperparameters, such as assigning a unique regularization weight to every individual parameter or learning rate schedules, which is computationally intractable for search-based methods. This advancement has materialized into practical tools, with open-source libraries such as \textit{Betty} \cite{choe2023betty}, \textit{Higher} \cite{grefenstette2019generalized}, \textit{TorchOpt} \cite{ren2023torchopt}, and \textit{JAXopt} \cite{blondel2022efficient} which leverage modern deep learning libraries that support automatic differentiation to enable efficient gradient-based hyperparameter optimization and meta-learning.


\subsection{Data poisining attacks}
For data poisoning attacks, it is important to note that in safety-critical scenarios, the leader (attacker) aims to maximize the learner's error by corrupting a subset of the training data. This creates a zero-sum game where the leader must guard against the follower (learner) selecting the robust solution $y \in S(x)$ that minimizes the damage. This corresponds to the pessimistic (or max-min) formulation:
\begin{equation} \label{eq:Pessimistic-Concrete}
    \max_{x \in X} \inf_{y \in S(x)} F(x, y) \coloneqq \sum_{\xi \in \mathcal{D}_{\text{val}}} \ell(y; \xi).
\end{equation}
In this case, the optimal solution set-valued mapping, $S: X \rightrightarrows \mathbb{R}^m$, of the training (lower-level) problem is defined by
\begin{equation}
     S(x) \coloneqq \argmin_{z \in \mathbb{R}^m} f(x, z) \coloneqq \sum_{\zeta \in \mathcal{D}_{\text{clean}}} \ell(z; \zeta) + \sum_{x_k \in x} \ell(z; x_k).
\end{equation}
Here, $x$ represents the set of poisoned data points injected by the attacker, and $x_k$ denotes the $k$-th individual poisoned sample within that set, whereas the lower-level variable $y$ represents the \textit{model parameters}. 
The set $X$ defines the feasible attack space, often bounded to ensure the poisoned data remains imperceptible or within valid input ranges. 

The pessimistic assumption ensures the attack is effective even against the best-case response of the learner. Seminal works in optimization-based poisoning include \cite{mei2015using, munoz2017towards}, while further pessimistic BL problems can be found in \cite{tian2021alphagan,lin2022evolutionary,liu2021value,huo2024perturbed,liu2023augmenting,
liu2021towards,ustun2024hyperparameter}. Gradient-based approaches for pessimistic BL are explored in \cite{guan2022gradient, huo2024perturbed}. 

It is worth noting that the literature often presents two sides of this coin. In \textit{Data Poisoning}, as formulated above, the Attacker is the Leader (maximizing error) and the Learner is the Follower. Conversely, in \textit{Adversarial Training}, the Learner is the Leader (minimizing error) who anticipates the Attacker (Follower) finding the worst-case perturbation. Both formulations are valid depending on whose perspective is being optimized. While the leader and follower clearly do not cooperate in these settings, many practical approaches rely on implicit function-based gradient methods. This is often the case because it is cheaper than using the pessimistic bilevel approach and it is mathematically justified when the lower-level problem has a unique solution (e.g., due to strong convexity or regularization). In such cases, the pessimistic and optimistic formulations coincide, $S(x)$ becomes a singleton, and the standard hypergradient derivations apply.
For more details on this subject, interested readers are referred to \cite{benfield2024classification,bruckner2011stackelberg,kantarciouglu2011classifier}.


\subsection{Further applications and structural challenges}
Beyond these examples, BO unifies diverse machine learning tasks \cite{liu2021investigating}. \textit{Meta-learning} fits this structure naturally, where the outer loop optimizes a meta-learner for fast adaptation on inner-loop tasks \cite{franceschi2018bilevel, finn2017model}.
Recently, BO has also become increasingly relevant for Large Language Models (LLMs). Architectures like \textit{Titans} \cite{behrouz2024titans} and the \textit{Nested Learning} framework \cite{behrouz2025nested} model long-term memory updates as an inner optimization loop. Similarly, \textit{Test-Time Training} \cite{sun2024learning, tandon2025end} enables adaptation to long contexts via autoregressive training on the input sequence during inference, and in this context  pretraining is similar to solving a BO problem where the lower-level is solved in an autoregressive fashion.
\textit{Deep equilibrium models} (DEQs) rely on fixed-point iterations, where training involves differentiating through the equilibrium state, a process mathematically equivalent to solving a bilevel program \cite{bai2019deep}. Similarly, \textit{generative adversarial networks} 
formulate a min-max game between a generator and discriminator, often treated as a bilevel program \cite{goodfellow2014generative, pfau2016connecting}.

The practical application of BO in machine learning faces significant structural difficulties. For the remainder of this section, we discuss some structural challenges associated to BL problems. We start with the \textit{high dimensionality and large-scale datasets.} A defining characteristic of modern machine learning is the scale of the problem. The lower-level variable $y$ is almost invariably high-dimensional. In deep learning, $y$ represents millions of neural network weights; in kernel methods, the number of dual variables scales with the dataset size. This makes it necessary to use approximate methods to compute the solution of the lower-level problem. The dimensionality of the upper-level variable $x$ varies: it is typically small in classical HPO ($<10$) but can be extremely large in Meta-learning or DEQs where $x$ represents network initializations or parameters.  Furthermore, the objective functions are defined as averages over massive datasets, precluding exact gradient evaluation and requiring the use of stochastic optimization methods. 

\paragraph{Non-uniqueness of the lower-level solution}
A central challenge is that the lower-level optimal solution set $S(x)$ is rarely a singleton for $x\in X$. The nature of this non-uniqueness varies by problem class. In the simplest case, such as linear regression with $\ell_2$ regularization (Ridge Regression), the lower-level objective $f$ is strongly convex with respect to $y$, guaranteeing a unique solution $S(x) = \{y(x)\}$. However, in the regime of \textit{over-parameterized} linear models—where the number of parameters $d$ exceeds the number of training samples $m_1$—if explicit regularization is absent, the objective is convex but not strongly convex. Consequently, $S(x)$ becomes an affine subspace containing infinitely many solutions that all achieve zero training error.
The situation becomes even more complex in \textit{deep learning}, where $f(x, y)$ is highly non-convex with respect to the neural network weights $y$ (see, e.g.,  \cite{li2018visualizing} for some graphical illustrations of the loss functions); 
In deep learning, non-uniqueness arises not just from flat valleys in the landscape but from fundamental symmetries (e.g., neuron permutation invariance) and the existence of multiple disconnected local minima. In this regime, the optimal solution set $S(x)$ is a disjoint union of manifolds, rendering the bilevel problem significantly harder to analyze.

\paragraph{The selection problem and algorithmic bilevel}
This non-uniqueness leads directly to the selection problem: when $S(x)$ contains multiple solutions, which one does the training algorithm actually return? Classical formulations assume the leader can select the best (optimistic) or worst (pessimistic) $y \in S(x)$. In practice, however, the solution $y$ is determined by the specific iterative algorithm used (e.g., Stochastic Gradient Descent) and its initialization. For instance, in over-parameterized linear regression initialized at zero, SGD converges specifically to the minimum Euclidean norm solution among the infinitely many minimizers. In deep learning, the ``implicit bias'' of the optimizer selects a specific attractor in the landscape. This algorithmic reality has led to ``unrolled'' formulations where the lower-level minimization is replaced by a dynamical system $y_{t+1} = \Phi_t(y_t, x)$ after $T$ steps, denoted $y_T(x)$. This perspective allows optimizing optimization hyperparameters (e.g., learning rates) and aligns the theoretical formulation with the algorithmic output \cite{grazzi2020iteration, arbel2022non, maclaurin2015gradient, franceschi2017forward}.

\paragraph{Generalization and the test set}
It is crucial to distinguish between the objective functions used in the bilevel optimization problem and the ultimate goal of the learning process. The upper-level objective $F$ is typically evaluated on a validation set $\mathcal{D}_{\text{val}}$ to guide the selection of $x$. However, this is merely a proxy; the true performance metric is the generalization error on a strictly held-out \textit{test set} $\mathcal{D}_{\text{test}}$, which is disjoint from both $\mathcal{D}_{\text{tr}}$ and $\mathcal{D}_{\text{val}}$ and is never used during the optimization of $x$ or $y$. This distinction highlights the risk of ``meta-overfitting'', where hyperparameters are tuned to minimize validation error but fail to generalize to the test distribution.

\paragraph{Constrained/nonsmooth/nonconvex lower-level problems}
Note that the HPO and data poisoning problems introduced above do not explicitly have lower-level constraints.  However, explicit lower-level constraints are fundamental to the formulation of hyperparameter tuning for Support Vector Machines \cite{bennett2006model, kunapuli2008classification}, fairness-aware learning \cite{yazdani2024fair}, sparse structure learning \cite{bertrand2020implicit, freconbilevel}, and safety-critical control \cite{wang2023bierl, kim2024spectral}. Further BL problems with lower-level constraints can be found in \cite{ye2023difference,gao2023moreau,liu2023value,yao2024constrained,xu2023efficient,xiao2023alternating,xiao2023alternatingOK,xiao2023alternatingPS,khanduri2023linearly,shi2023double,liu2021towards,liu2021value,liu2022bome,ye2024first}. These constraints necessitate the handling of non-differentiable projection operators and the breakdown of strict complementarity in KKT conditions \cite{bolte2022nonsmooth}, which complicate the estimation of hypergradients. Recent algorithmic literature has addressed these issues in various ways. For example using one-stage primal-dual formulations that update variables simultaneously \cite{jiang2024primal, sow2022primal}, smoothed implicit gradient approximations \cite{khanduri2023linearly}, and difference-of-convex relaxations for value-function constraints \cite{gao2022value}. 
In BL, the lower-level problem can be nonsmooth (see, e.g., \cite{liu2023rethinking,shi2022efficient,okuno2021,alcantara2021unified,bertrand2020implicit,bertrand2022implicit}) or even nonconvex (see, e.g., \cite{arbel2022non,liu2023value,liu2021value,liu2021towards,liu2023augmenting,liu2022bome,lu2024slm}). In such scenarios, the classical algorithmic framework is not applicable, as it will be clear in Sections \ref{sec:Main algorithmic techniques} and \ref{sec:Challenges and limitations}.

\section{Implicit function-based methods for bilevel learning}\label{sec:Main algorithmic techniques}
The implicit function model \eqref{eq:Implicit-Function-Model} is the framework for the classical algorithms for BL in the current literature. In this section, we examine the main techniques that have been developed so far from this perspective. 
More precisely, the main algorithmic idea is the gradient descent method tailored to the specific version 
\begin{equation}\label{eq:bilevellls}
\min_{x \in \mathbb{R}^n}~\mathcal{F}(x):=F(x, y(x)) \quad \mbox{ with }\quad  y(x) := \arg \min_{y \in \mathbb{R}^m}~f(x,y)
\end{equation}
of problem \eqref{eq:Implicit-Function-Model}, where in addition to the requirement that the lower-level problem has a unique optimal solution (for all upper-level variables), it is  assumed  that there is no upper-- nor lower--level constraints (i.e., with $X:=\mathbb{R}^n$ and $Y(x):=\mathbb{R}^m$). While upper-level constraints are generally manageable, incorporating lower-level constraints presents additional challenges which will be discussed in later sections. This formulation is well-defined under the basic assumption that $y(\cdot): \mathbb{R}^n\rightarrow \mathbb{R}^m$ describes a vector-valued function such that $\{y(x)\}=S(x)$ (the optimal solution set is a singleton) for all $x\in \mathbb{R}^n$.

To have a sense of how this is still a very challenging problem, start by noting that we can easily  find examples where the upper-level objective function $F$ is convex and the lower-level objective function is strongly convex in the lower-level variable $y$, but the resulting function $x\mapsto \mathcal{F}(x)$ is still not necessarily convex. For example, let 
\begin{equation*}\label{eq:Example-1}
    F(x,y) := y \;\; \mbox{ and } \;\; f(x,y) := \left(y-(x^3-x)\right)^2.
\end{equation*}
Clearly, $f(x, \cdot)$ is strongly convex for all $x\in \mathbb{R}$, and we can observe that the resulting function $x \mapsto\mathcal{F}(x):= F(x,y(x))=x^3 - x$ is nonconvex. Hence, \eqref{eq:Implicit-Function-Model} is a nonconvex optimization problem, even when all necessary requirements are satisfied for it to be well-defined. Additionally, as it will be clear in Section \ref{sec:Constrained optimization}, the assumptions required for \eqref{eq:Implicit-Function-Model} to be well-defined are generally only local. This means in many practical use cases the model is valid only locally and not necessarily on the whole domain of the reduced objective function $\mathcal{F}$ or the upper-level feasible set $X$.

We can also consider the stochastic setting where the upper- and lower-lever objective functions are expectations, which model situations where we have access to data sampled from unknown distributions. For such a setup, we have
\begin{equation}\label{eq:bilevelllsstoch}
F(x, y) =\mathbb{E}_\xi \hat F(x, y, \xi) \quad \mbox{ and } \quad f(x, y) = \mathbb{E}_\zeta \hat f(x,y, \zeta),
\end{equation}
where $\xi$ and $\zeta$ are random variables. Some works study also the special case where expectations are replaced with their empirical estimates, i.e., they are finite sums 
\begin{equation}\label{eq:bilevelllsfinite}
F(x, y) =\frac{1}{d_u}\sum_{i=1}^{d_u} \hat F(x, y, \xi_i) \quad \mbox{ and } \quad  f(x, y) = \frac{1}{d_l}\sum_{i=1}^{d_l} \hat f(x,y, \zeta_i).
\end{equation}
As we mentioned in \Cref{sec:A flavor of machine learning}, in BL, $n$, $m$, $d_u$ and $d_l$ can be very large and $f$ and $F$ might be nonconvex: the main example is whenever the lower-level and/or upper-level variables are the parameters of a neural network.  In the case where the function $f$ is nonconvex w.r.t. $y$, the optimal solution set of the lower-level problem is not necessarily a singleton. However, algorithms are still applied in this setup with various degree of success.

\subsection{General algorithmic framework}
The core idea of the implicit function approaches is to compute the \textit{hypergradient}, i.e.\@ the gradient of $\mathcal{F}$ by exploiting the Implicit Function Theorem, which characterizes the dependence of the lower-level optimal solution on the upper-level variables.
More precisely, if for a given $x$ the lower-level problem admits a \emph{unique} local minimizer $y(x)$, and if $f$ is twice continuously differentiable with $\nabla^2_{yy} f(x,y(x))$ invertible, then the first-order optimality condition
\begin{equation}\label{eq:nabla_y_f(x,y)}
   \nabla_y f(x, y(x))=0, 
\end{equation}
defines $y$ locally as an implicit function of $x$. Differentiating this stationarity condition with respect to $x$ gives
\[
\nabla^2_{xy} f(x,y(x))+\nabla^2_{yy} f(x,y(x))\,\nabla y(x)=0,
\]
and, by the Implicit Function Theorem,
\begin{align*}    
\nabla y(x)&= -\left[\nabla^2_{yy}f(x,y(x))\right]^{-1}\nabla^2_{xy}f(x,y(x)), \\
\nabla \mathcal{F}(x)
&=
\nabla_x F(x,y(x))
+
\nabla y(x)^\top \nabla_y F(x,y(x)).
\end{align*}

The gradient-based method that computes the exact hypergradient is described in Algorithm \ref{alg:1}; see, e.g., \cite{dempe2002foundations}.

\begin{algorithm}
\caption{Classical hypergradient descent algorithm}
\label{alg:1}
\begin{algorithmic}[1]
\REQUIRE $x^0$ and sequence $\{\alpha_k\}$.
\FOR{$k=0, 1, \dots, K$}
\STATE  $y(x^k)=\arg\underset{y}{\min}f(x^k, y)$
\STATE  $\nabla y(x^k) = -\left[\nabla^2_{yy}f(x^k, y(x^k))\right]^{-1} \nabla^2_{xy}f(x^k, y(x^k))$
\STATE  $\nabla\mathcal{F}(x^k)=\nabla_x F(x^k, y(x^k)) + \nabla y(x^k)^\top \nabla_y F(x^k, y(x^k))$
\STATE $x^{k+1}=x^k - \alpha_k \nabla\mathcal{F}(x^k)$
\ENDFOR
\end{algorithmic}
\end{algorithm}

The method relies on the well-defined nature of problem \eqref{eq:Implicit-Function-Model}, where the Jacobian $\nabla y(x^k)$ is computed by inverting the Hessian of the lower-level objective function. The commonly used sufficient condition to make the Jacobian well-defined is for the lower-level objective function to be strongly convex and twice continuously differentiable, which is satisfied, for example, by having an $\ell_2$ regularization penalty in combination with a convex and twice differentiable objective such as in regularized logistic regression. 

In high-dimensional settings, however, Lines 2 and 3 of \Cref{alg:1} are computationally prohibitive. First, computing the lower-level solution $y(x^k)$ typically requires running an iterative lower--level solver to (near) convergence; in large-scale problems this can mean many gradient/Hessian--vector evaluations (often over large datasets), so the cost scales with both the lower-level dimension and the number of lower--level iterations. Second, the implicit-gradient term involves inverting the hessian. Forming the hessian already costs $\mathcal{O}(m^2)$ memory/time, and a direct factorization (e.g., Cholesky/LU) costs $\mathcal{O}(m^3)$ time and $\mathcal{O}(m^2)$ memory, which quickly becomes infeasible as $m$ grows.

It is important to acknowledge that alternative frameworks exist which bypass the direct computation of the hypergradient. Notably, more recent efficient methods like BOME \cite{liu2022bome} adopt a value-function approach (See \Cref{sec:Pure first and second order methods}), reformulating the bilevel program into a single-level constrained optimization problem. This allows for the application of fully first-order methods that do not require computing second derivatives of the lower--level objective.

\subsection{Efficient approximation schemes}
To address the computational bottlenecks of Algorithm \ref{alg:1},  various approximation schemes have been introduced. These methods generally aim to estimate the hypergradient efficiently without incurring the full cost of Hessian inversion or exact lower-level minimization.

A fundamental distinction in the literature is between \textit{Approximate Implicit Differentiation (AID)} and \textit{Iterative Differentiation (ITD)}. While AID approximately solves the linear system derived from the expression of the hypergradient, 
ITD approximates the hypergradient by backpropagating through the unrolled computational graph of the lower-level optimization algorithm (e.g., $T$ steps of gradient descent). The two procedures are detailed in \Cref{alg:aid,alg:itd}, respectively. As explained in \cite{maclaurin2015gradient,franceschi2018bilevel}, ITD computes the exact gradient of the proxy objective $x \mapsto F(x, y_T(x))$. However, its memory cost scales linearly with $T$, which can be prohibitive for problems requiring many lower--level steps. In contrast, AID allows for constant memory cost (using efficient linear solvers) but introduces a systematic bias if the linear system is not solved exactly or if the lower-level problem has not converged to the true solution $y^*(x)$. For a detailed comparison between these two strategies, interested readers are referred to \cite{grazzi2020iteration}.

\begin{algorithm}
    \caption{Bilevel Gradient Method with AID}
    \label{alg:aid}
    \begin{algorithmic}[1]
    \REQUIRE $x^0$, $\{\alpha_k\}$.
    \FOR{$k=0, 1, \dots, K$}
        \item[] \textbf{Lower--level optimization:}
        \STATE Init solver at $\hat{y}_0$ (if warm start $\hat{y}_0 = \hat{y}_{k-1}$)
        \STATE Find $\hat{y}_k \approx \arg\min_y f(x^k, y)$
        \item[] \textbf{Linear system:}
        \STATE Init solver at $\hat{v}_0$ (if warm start $\hat{v}_0 = \hat{v}_{k-1}$)
        \STATE Find $\hat{v}_k$ by solving  $\nabla^2_{yy}f(x^k, \hat{y}_k) v = \nabla_y F(x^k, \hat{y}_k)$ 
        \item[] \textbf{Approximate Hypergradient:}
        \STATE $h_k = \nabla_x F(x^k, \hat{y}_k) - \nabla^2_{xy}f(x^k, \hat{y}_k)^\top \hat{v}_k$
        \item[] \textbf{Upper--level update:}
        \STATE $x^{k+1} = x^k - \alpha_k h_k$
    \ENDFOR
    \end{algorithmic}
\end{algorithm}

In the AID framework, the computation of the hypergradient approximation (see Lines 4 and 5 in Algorithm \ref{alg:aid}) is typically reformulated using an adjoint vector $v$, which is the solution to the linear system $\nabla^2_{yy}f(x^k, \hat y_k) v = \nabla_y F(x^k, \hat y_k)$, so that 
\[
\nabla\mathcal{F}(x^k) \approx h_k =\nabla_x F(x^k, \hat y_k) + \nabla^2_{xy}f(x^k, \hat y_k)^\top v.
\]
This linear system can be solved for example using Conjugate Gradient (CG) \cite{pedregosa2016hyperparameter}, which converges fast and relies solely on Hessian-vector products and avoids explicit matrix factorization. Alternatively, the linear system can be viewed as the optimality condition of a quadratic minimization problem and solved using (stochastic) gradient descent \cite{grazzi2020iteration, grazzi2021convergence}. This generalizes the Neumann series approximation (of the inverse hessian) \cite{lorraine2020optimizing}, which corresponds to a specific gradient descent scheme. 

To further reduce computational overhead, recent ``Hessian--free'' methods attempt to mitigate the cost of accessing second-order information directly. Unlike methods that simply ignore second-order terms, Hessian-free bilevel algorithms approximate the required Hessian-vector products to maintain convergence to the true stationary point. As analyzed by Sow et al. \cite{sow2022hessian}, one can estimate the product $\nabla^2_{yy} f(x, y) v$ using a finite difference approximation of gradients for a sufficiently small scalar $\delta > 0$. This technique allows the linear system to be solved using only first-order gradient oracles, avoiding explicit Hessian construction while strictly adhering to the implicit differentiation framework. 

\begin{algorithm}
    \caption{Bilevel Gradient Method with ITD}
    \label{alg:itd}
    \begin{algorithmic}[1]
    \REQUIRE $x^0$, $\{\alpha_k\}$, steps $T$.
    \FOR{$k=0, 1, \dots, K$}
        \item[] \textbf{Lower--level optimization:}
        \STATE Init $y_0$ (if warm start set $y_{0}(x_k) = y_{T}(x_{k-1})$)
        \FOR{$t=0, \dots, T-1$}
            \STATE $y_{t+1}(x^k) = \Phi_t(x^k, y_t(x_{k}))$ \COMMENT{E.g., $\Phi_t(x, y) = y-\eta_t\nabla_y f(x,y)$}
        \ENDFOR
        \STATE \text{Set } $\mathcal{F}_T(x^k) := F(x^k, y_T(x^k))$
        \item[] \textbf{Approximate Hypergradient:}
        \STATE \text{Compute } $h_k = \nabla \mathcal{F}_T(x^k)$ \text{via backpropagation through $\Phi_0, \dots, \Phi_{T-1}$}
        \item[] \textbf{Upper-level update:}
        \STATE $x^{k+1} = x^k - \alpha_k h_k$
    \ENDFOR
    \end{algorithmic}
\end{algorithm}

In the very common situation where we are dealing with large scale datasets, i.e.\@ where $d_u, d_l$ in \eqref{eq:bilevelllsfinite} are large, the lower-level, the linear system and the final hypergradient computation (for AID methods),  are solved with iterative algorithms relying on the \textit{stochastic estimators}  $\hat f$ and $\hat F$ of the lower and upper level objectives and their derivatives.

A ubiquitous strategy to improve efficiency is ``warm starting'', which involves initializing the lower-level and linear system solvers at iteration $k$ with the approximate solutions from iteration $k-1$. Intuitively, this exploits the smoothness of the solution map $y^*(x)$, implying that small changes in $x$ result in small changes in the optimal $y$, making the previous solution a high-quality initial guess. Recent theoretical work by \cite{ji2022will} has rigorously analyzed the impact of the computational ``loops'' (lower--level optimization steps $T$ and linear solver steps $Q$) on convergence. They demonstrate that while substantial loops are often necessary for ITD, AID schemes utilizing warm-start strategies can effectively reduce the required loop count to $O(1)$ per iteration. The special case where only one step is done for the lower-level problem is referred to as single-loop (see e.g. \cite{li2022fully,gong2024nearly}), in contrast with double-loop bilevel algorithms which use more than one step. While this significantly improves practical efficiency, it complicates the theoretical analysis by introducing a strong coupling between the dynamics of the upper and lower-levels.

\subsection{Convergence guarantees and complexity analysis}

Formalizing the convergence of bilevel algorithms requires distinguishing between the quality of the \textit{hypergradient approximation} and the convergence of the \textit{optimization procedure} as a whole.

\subsubsection{Convergence of the hypergradient approximation}

Before analyzing the optimization trajectory, one must ensure that the estimated hypergradient $h_k$ is a reliable proxy for the true hypergradient $\nabla \mathcal{F}(x_k)$. The true hypergradient depends on the exact lower-level solution $y^*(x_k)$ and the exact solution to the linear system $v^*(x_k)$, while in practice, we operate with approximations.


\paragraph{Lower-level fixed points}
\cite{grazzi2020iteration,grazzi2021convergence} study the accuracy of (approximate) hypergradient computation for both AID and ITD bilevel problems in which the lower-level solution is defined implicitly as a fixed point,
$y(x)=\Phi(x,y(x))$, under regularity assumptions that ensure well-posedness and stability of the implicit map.
This framework subsumes the classical smooth lower level optimality condition $\nabla_y f(x,y(x))=0$ by choosing $\Phi(x,y)$ as a single step of a descent method, e.g., $\Phi(x,y)=y-\eta \nabla_y f(x,y)$, so that fixed points coincide with stationary points of the lower-level objective.

Importantly, the fixed-point viewpoint is strictly more general than convex optimization in the sense of minimizing a scalar potential: many equilibrium problems are naturally modeled by \emph{monotone operator} formulations and solved by \emph{operator-splitting} iterations (proximal point, forward--backward, Douglas--Rachford, etc.), which can be written as fixed-point iterations even when the underlying operator does \emph{not} arise as the gradient (or subgradient) of any function \cite{ryu2016primer,bauschke2020correction}.
Intuitively, an operator fails to come from an optimization problem when it lacks a potential structure (e.g., it contains an intrinsically ``rotational''/skew-symmetric component), as happens in general variational inequalities, game-theoretic equilibria, and saddle-point problems where the relevant stationarity conditions correspond to a monotone inclusion rather than minimization of a single scalar objective.
This also connects to \emph{implicit} or \emph{equilibrium} architectures in deep learning, notably Deep Equilibrium Models \cite{bai2019deep}, which define hidden states as solutions of a fixed-point equation.

Under the assumption that  $\Phi$ (or $\nabla f$) and $F$ are Lipschits smooth, and that $\Phi$ is a \textit{contraction} (Lipschitz with constant less than one) with respect to $y$ (or alternatively if $f$ is smooth and strongly convex with respect to $y$), the error in the AID hypergradient is bounded linearly by the errors in the lower-level and linear system subproblems:
\begin{equation} \label{eq:grazzi_bound}
    \| h_k- \nabla \mathcal{F}(x_k) \| \leq L_1 \| \hat{y}_k - y^*(x_k) \| + L_2 \|  \hat{v}_k - v^*(x_k) \|,
\end{equation}
where $L_1, L_2$ are constants derived from the condition number $\kappa$ of the lower-level problem and the Lipschitz constants of the gradient and/or fixed-point maps. Specifically, if 
we employ gradient descent for the lower-level and Conjugate Gradient (CG) or gradient descent for the linear system, \cite{grazzi2020iteration} show that achieving an error $\| \nabla \mathcal{F}(x_k) - h_k \| \leq \epsilon$ requires $O(\log(1/\epsilon))$ iterations for both the lower-level solver and the linear system solver.

Under the same contractivity and smoothness assumptions, \cite{grazzi2020iteration} also establish a \emph{linear} (geometric) convergence of the ITD hypergradient to the true hypergradient. A similar results can be also found in the seminal automatic differentiation textbook by Griewank and Walther \cite[Chapter 15]{griewank2008evaluating}. In particular, for $T$ lower--level gradient descent steps (unrolling length), \cite[Theorem~2.1]{grazzi2020iteration} shows a bound of the form
\[
\bigl\|\nabla \mathcal{F}_T(x^k)-\nabla \mathcal{F}(x^k)\bigr\|
\;\le\;
(c_1T +c_2)\,q^{\,T},
\]
where $c_1,c_2 \geq 0$ and $0\leq q< 1$ depend on $x$. This rate similar to the AID rate of $cq^T$, but multiplied by an additional prefactor growing like $T$.
As a consequence, the iteration complexity to reach accuracy $\epsilon$ remains $O(\log(1/\epsilon))$, but the extra $T$ prefactor yields a practical ``delay'' compared to AID bounds.


The \textit{stochastic case} has been mainly studied for AID. The hypergradient estimator $h_k$ admits an MSE decomposition into (i) terms controlled by the mean-square accuracy of the stochastic lower--level solvers used to compute $(\hat y_k,\hat v_k)$, and (ii) a variance term coming from stochastic Jacobian/Hessian-vector product estimators \cite{grazzi2021convergence,grazzi2023bilevel}.
As a result, the overall MSE essentially matches the convergence rate of the stochastic solvers employed at the lower-level and for the linear system, up to an additional variance contribution. In particular, the refined analysis in \cite{grazzi2023bilevel} yields bounds of the form
\[
\mathbb{E}\|h_k-\nabla\mathcal{F}(x_T)\|^2 \;=\; O\!\left(\rho_T + \sigma_T + \frac{1}{b_T}\right),
\]
where $\rho_T$ and $\sigma_T$ denote the MSEs of the lower-level and linear-system stochastic solvers (respectively), and $b_T$ is number of samples used for some stochastic estimators. Hence, when both lower--level solvers are implemented with SGD-type methods so that $\rho_T=\Theta(1/T)$ and $\sigma_T=\Theta(1/T)$ and we also set $b_T = \Theta(T)$, the hypergradient MSE decreases as $O(1/T)$ as well, matching the canonical $O(1/T)$ rate of SGD for the lower-level problem. The ITD case has been studied only  recently by \cite{iutzeler2024derivatives}, where they focus on the convergence of the derivative of SGD. They show that with a careful analysis it is possible to achieve a rate of $O(\log(T)^2/T)$, which is off only by a logarithmic factor compared to the AID rate.

While the classical theory relies on differentiability, recent work has extended the implicit differentiation framework to nonsmooth settings, such as when the lower-level problem involves $L_1$ regularization or nonsmooth activations. \cite{bolte2022nonsmooth} and \cite{bertrand2022implicit} established that under mild conditions (e.g., separability or specific non-degeneracy), implicit differentiation remains valid on the support of the solution, and convergence rates can still be derived. \cite{grazzi2024nonsmooth} further generalized the analysis to provide non-asymptotic convergence rates for nonsmooth  AID, including the stochastic case.

\subsubsection{Convergence of the bilevel framework}

\begin{table*}[htbp]
    \centering
    \setlength{\tabcolsep}{4pt}
    \renewcommand{\arraystretch}{1.2}
    \resizebox{\textwidth}{!}{
        \begin{tabular}[t]{clcc}
            \toprule
            \multicolumn{4}{c}{\textbf{Stochastic Setting}} \\
            \midrule
            \textbf{Complexity} & \textbf{Algorithm} & \textbf{WS} & \textbf{Loop} \\
            \midrule
            $O(\epsilon^{-3})$ & BSA \cite{ghadimi2018approximation} & No & Nested \\
            \midrule
            $\tilde{O}(\epsilon^{-2.5})$ & TTSA \cite{hong2023two} & Yes & Single \\
            \midrule
            \multirow{3}{*}{$\tilde{O}(\epsilon^{-2})$} 
            & stocBiO \cite{ji2021bilevel} & Yes & Nested \\
            & SMB \cite{guo2021randomized}, saBiAdam \cite{huang2021biadam} & Yes & Single \\
            & ALSET \cite{chen2021tighter} & Yes & Single \\
            \midrule
            \multirow{2}{*}{$O(\epsilon^{-2})$} 
            & BSGM \cite{grazzi2023bilevel} & No & Nested \\
            & Amigo \cite{arbel2022amortized} & Yes & Nested  \\
            \specialrule{1.2pt}{3pt}{3pt}
            \multicolumn{4}{c}{\textbf{Variance Reduction}} \\
            \midrule
            {$O(\epsilon^{-2})$}  & STABLE \cite{chen2021single}, FSLA \cite{li2021fully} & Yes & Single \\
            \midrule
            \multirow{3}{*}{$\tilde{O}(\epsilon^{-1.5})$} 
            & VRBO \cite{yang2021provably} & Yes & Nested \\
            & STABLE-VR \cite{guo2021randomized}, SUSTAIN \cite{khanduri2021near}, & \multirow{2}{*}{Yes} & \multirow{2}{*}{Single} \\
            & VR-saBiAdam \cite{huang2021biadam}, MRBO \cite{yang2021provably} & & \\
            \midrule
            $O(d^{2/3}/\epsilon)$ & SABA \cite{dagreou2024lower} & Yes & Nested \\
            \bottomrule
        \end{tabular}%
        \hspace{0.02\textwidth}%
        \begin{tabular}[t]{clcc}
            \toprule
            \multicolumn{4}{c}{\textbf{Deterministic Setting}} \\
            \midrule
            \textbf{Complexity} & \textbf{Algorithm} & \textbf{WS} & \textbf{Loop} \\
            \midrule
            $O(\epsilon^{-5/4})$ & BA \cite{ghadimi2020approximation} & No & Nested \\
            \midrule
            \multirow{2}{*}{$\tilde O(\epsilon^{-1})$} 
            & BiO-ITD \cite{ji2021bilevel} & Yes & Nested \\
            & BGM \cite{grazzi2023bilevel} & No & Nested \\
            \midrule
            \multirow{2}{*}{$O(\epsilon^{-1})$} 
            & BiO-AID \cite{ji2021bilevel} & Yes & Nested \\
            & Amigo \cite{arbel2022amortized} & Yes & Nested \\
            \bottomrule
        \end{tabular}
    }
    
    \vspace{0.3cm}
    
    \caption{Sample complexity for finding an $\epsilon$-stationary point (i.e. a point $x$ such that $\mathbb{E} \| \nabla \mathcal{F}(x) \|^2 \leq \epsilon$) in implicit function based methods. \textit{WS} indicates the use of warm-start for the lower-level problem. \textit{Loop} specifies the structure: ``Single'' implies strictly 1 lower--level step is taken per upper--level step, while ``Nested'' implies a multi-step lower--level solver is used. The complexity of SABA is for the finite sum setting where $d=d_u + d_l$ is the total number of samples.}
    \label{tab:implicit_function_complexity}
\end{table*}

To analyze the complexity of the full framework, we must first define the computational oracles involved. The analysis typically revolves around following operations: the \textit{lower and upper level gradients} $\nabla_y f(x, y)$, $\nabla_x \mathcal{F}(x,y)$ ,$\nabla_y \mathcal{F}(x,y)$ the \textit{lower-level Hessian-vector product} $\nabla^2_{yy} f(x, y) v$, and the \textit{mixed Hessian-vector product} $\nabla^2_{xy} f(x, y)^\top v$. In the stochastic case, these quantities are replaced instead by unbiased estimates with controlled variance. Thanks to Automatic Differentiation, computing these derivative operations has a cost which scales with the dimension of $x$ and $y$ in the same way as computing the function values. While the Hessian-vector product is typically implemented as a Jacobian-Vector Product (JVP) of the lower-level gradient map with respect to $y$, practical implementations often find that simple gradient evaluations are significantly cheaper (by a constant factor) than the second-order vector products. For convergence rates, these costs are often aggregated, measuring the total number of calls to these first-order and second-order oracles required to reach stationarity.

Since the the function $\mathcal{F}$ is nonconvex, the goal of the bilevel procedure is to find a stationary point. Therefore, the standard convergence metric is the number of iterations (or total oracle calls) required to find an $\epsilon$-stationary point, defined as a point $x$ satisfying $\mathbb{E} \| \nabla \mathcal{F}(x) \|^2 \leq \epsilon$. We summarize the main results in \Cref{tab:implicit_function_complexity} and discuss them in detail below.
In the deterministic case, the seminal work by \cite{ghadimi2020approximation} established a baseline complexity rate of $O(\epsilon^{-5/4})$ without warm start for AID. Subsequent research improved upon this by leveraging the warm-start strategy. By initializing the solvers with previous iterates, \cite{ji2021bilevel} and \cite{arbel2022amortized} demonstrated that it is possible to achieve the optimal complexity of $O(\epsilon^{-1})$, which matches the standard rate for single-level nonconvex optimization. The work \cite{ji2021bilevel} also establishes the first rate for ITD of $O(\epsilon^{-1}\log(\epsilon^{-1}))$ for ITD, which relies on warm-start.

The stochastic setting is significantly more challenging due to the bias-variance trade-off and the results are focused almost entirely on the AID method. The landscape here is defined by the sample complexity required to reach an $\epsilon$-stationary point. The baseline complexity was established by the Bilevel Stochastic Approximation (BSA) algorithm \cite{ghadimi2020approximation}, which does not use warm start and requires $O(\epsilon^{-3})$ samples and an increasing (as $\sqrt{k}$ where $k$ is the upper-level iteration counter) number of lower--level steps to control bias. Subsequent improvements utilized warm start to reduce the lower--level loop complexity. The Two-Timescale Stochastic Approximation (TTSA) \cite{hong2023two} achieves $\tilde{O}(\epsilon^{-2.5})$. 

A significant leap in efficiency was marked by a class of algorithms achieving $\tilde{O}(\epsilon^{-2})$ complexity, matching the standard single-level stochastic baseline. This group includes stocBiO \cite{ji2021bilevel}, SMB \cite{guo2021randomized}, ALSET \cite{chen2021closing}, Amigo \cite{arbel2022amortized}, STABLE \cite{chen2022stable}, and FSLA \cite{li2022fully}. Most of these methods achieve this efficiency by using warm start on the lower-level problem to reduce the number of lower--level iterations to $O(1)$ (often referred to as single-loop), although some, like Amigo, warm-start both the lower-level and linear system solvers which allows to remove the log factor and achieve $O(\epsilon^{-2})$ sample complexity, which interestingly matches that of a single level optimization problem. Variance-reduced methods such as SUSTAIN \cite{khanduri2021near}, VR-saBiAdam \cite{huang2021biadam}, and MRBO \cite{yang2021provably} push the theoretical boundary further, achieving a near-optimal complexity of $\tilde{O}(\epsilon^{-1.5})$, though they typically require more oracles per iteration and stronger assumptions.
Specifically focusing on the finite-sum setting (Empirical Risk Minimization), \cite{dagreou2024lower} established tight lower bounds and proposed the SABA algorithm, which achieves optimal variance-reduced rates. This approach builds on their earlier general stochastic framework \cite{dagreou2022framework}, which provided a unified analysis for single-loop aggregation schemes applicable to both finite-sum and infinite-horizon problems.

While warm start is prevalent, \cite{grazzi2023bilevel} challenged the assumption that it is necessary for optimal rates. They showed that a ``cold start'' procedure (solving the lower--level problem from scratch or fixed initialization) can still achieve the $\tilde{O}(\epsilon^{-2})$ sample complexity in the stochastic setting. The key is using larger mini-batches ($\Theta(\epsilon^{-1})$) for the hypergradient estimation and specific step-size schedules to control the variance. For the deterministic case, this cold-start approach improves the baseline to $O(\epsilon^{-1}\log(\epsilon^{-1}))$. While without warm-start the analysis is simpler since it decouples the lower--level and outer loop dynamics, it additionally requires the distance between the lower-level solution and the starting point of the lower-level solver to be uniformly bounded over the upper-level feasible set, which might explain why warm-start method often perform better in practice.

\subsection{Meta-parameters and adaptive methods}\label{sec:adaptive}
It is worth emphasizing that many bilevel algorithms with strong convergence guarantees depend on several \emph{algorithmic meta-parameters}, that is, quantities controlling how the numerical method is run rather than parameters of the learning model itself. For example the upper- and lower-level stepsizes, the number of lower-level updates performed per outer iteration, mini-batch sizes in the stochastic setting, regularization or damping parameters in implicit solvers, and tolerances or number of iterations for auxiliary linear-system solves. In many methods, these quantities cannot be chosen once and kept fixed, but instead must follow carefully designed schedules along the iterations in order to balance bias, variance, and stability. This tuning burden has motivated a recent line of work on adaptive and parameter-free bilevel methods. Early work in this direction includes \emph{BiAdam}, which introduced adaptive learning rates for stochastic bilevel optimization and its variance-reduced variant VR-BiAdam \cite{huang2021biadam}. 

More recently, \emph{BiSLS/SPS} proposed stochastic line-search and Polyak-type rules to automatically choose the coupled upper- and lower-level stepsizes, with the explicit goal of improving stability and reducing manual tuning \cite{fan2023bisls}. Going one step further, tuning-free methods such as D-TFBO and S-TFBO remove the need for prior knowledge of problem constants and choose stepsizes adaptively from cumulative gradient information, while attaining deterministic and stochastic rates that nearly match those of their well-tuned counterparts \cite{yang2024tuning}. Related adaptive mirror-descent variants have also been proposed beyond the strongly-convex inner-level setting \cite{antonakopoulos2025adaptive}. 
Adaptive ideas have also started to extend beyond the standard centralized Euclidean setting, including federated, Riemannian, and decentralized bilevel optimization, further indicating that robustness to solver meta-parameters is emerging as a broader theme in the bilevel literature \cite{huang2022federated,shi2025adaptive,zhai2025problem}.

\subsection{Beyond classical assumptions}
Some works have begun to relax classical assumptions in BL, to extend the gradient descent framework presented here to more general problem classes. In particular regarding nonconvex geometry, weak lower-level curvature, and fundamental complexity limits.
Classical results focus on finding first-order stationary points, 
characterized solely by a small gradient norm. However, in nonconvex BO landscapes, such points may correspond to strict saddles of the upper-level objective function. To address this,  \cite{huang2022efficiently} develops perturbed implicit-gradient methods that guarantee convergence to second-order stationary points. 
Common BO approaches to tackle problems with nonconvex lower-level players are discussed in Section \ref{sec:Pure first and second order methods}.

The \textit{lower-level singleton} 
assumption fails in many over-parameterized models, where the lower-level optimal solution may not be unique and the Hessian of the lower-level objective function is singular. To address this, recent papers have introduced Polyak-Łojasiewicz (PL) conditions for the lower-level problem. \cite{liu2021towards} propose a gradient-based framework that handles nonconvex followers by utilizing initialization auxiliaries, bypassing the need for strong convexity. Possible approaches to extend the gradient descent framework to problems where the lower-level singleton assumption (described in \eqref{eq:S(x)=1}) fails are discussed in the next section. 
Furthermore, \cite{xiao2023generalized} and \cite{huang2023momentum} extend convergence guarantees to settings satisfying the PL condition, showing that implicit differentiation strategies remain tractable even when the argmin set is non-singleton. These works often rely on generalized inverses (such as the Moore--Penrose pseudoinverse) or iterative approximations to handle the singular Hessian matrices inherent in these relaxed settings.

Recent lower-bound results demonstrate that BO is fundamentally harder than minimax or single-level nonconvex optimization. \cite{ji2023lower} establish oracle lower bounds that match (up to logarithmic factors) the best known upper bounds, showing the intrinsic dependence on conditioning and cross-level coupling. In parallel, \cite{dagreou2024lower} prove tight lower bounds for bilevel empirical risk minimization and introduce an algorithm achieving near-optimal complexity.

\section{Limitations of the implicit function-based framework}\label{sec:Challenges and limitations}

The gradient descent algorithmic framework discussed in the previous section has been very successful in solving various BL problem classes as highlighted previously. However, the assumptions for its implementation are very restrictive. At a  high-level, the main restrictions are the requirement that \textit{the lower-level player is unconstrained and can  only have a unique optimal solution for all upper-level variables}. In this section, we analyze possible ways to extend the gradient descent method to the implicit function-type framework, while relaxing these assumptions. We start by first keeping assumption \eqref{eq:S(x)=1}, while relaxing the assumption that the lower-level problem is unconstrained (i.e., the requirement that $Y(x) = \mathbb{R}^m$ for all $x\in X$).  

\subsection{Extension of the implicit function approach to lower-level constrained  problems} We assume here that the lower-level feasible set has the form
\begin{equation}\label{eq:Y(x)_g(x,y)}
    Y(x):=\left\{y\in \mathbb{R}^m\left|\;\, g(x, y)\leq 0\right.\right\},
\end{equation}
where $g: \mathbb{R}^n \times \mathbb{R}^m \rightarrow \mathbb{R}^q$  corresponds to  the lower-level constraint function. 
Recall that the key requirement in the previous section to ensure the fulfillment of assumption \eqref{eq:S(x)=1} is to assume  that the function $f(x, \cdot)$ is strongly convex for all $x\in X$. Next, we introduce assumptions that can ensure \textit{strong local stability} for the lower-level problem; i.e.,  to enable the reduction of $y(\cdot)$ to a locally Lipschitz continuous vector-valued and even differentiable function,  when the lower-level feasible is given by \eqref{eq:Y(x)_g(x,y)}. Before, note that in Section \ref{sec:Constrained optimization} (resp. Section \ref{sec:Pure first and second order methods}), we are going to specifically explore classical BO approaches to tackle bilevel programs where the lower-level problem is only convex (resp. nonconvex).  

To ensure that assumption \eqref{eq:S(x)=1} holds under the lower-level feasible set \eqref{eq:Y(x)_g(x,y)}, we can make the following assumptions associated to a  point $(\bar x, \bar y)$ such that $\bar y\in Y(\bar x)$:
\begin{itemize}
\item[(A1)] The function $f$ and $g$ are at least twice continuously differentiable near $(\bar x, \bar y)$.
\item[(A2)] The lower-level linear independence constraint qualification (LLICQ) holds at $(\bar x, \bar y)$; i.e., the vectors $\nabla_y g_i(\bar x, \bar y)$ for $i\in I(\bar x, \bar y)$ are linearly independent. Here, $I(\bar x, \bar y)$ denotes the set of active indices for the lower-level constraint:
    \[
    I(\bar x, \bar y):=\left\{i\in[q]\,|\;\; g_i(\bar x, \bar y)=0\right\} \; \mbox{ with } \; [q]:=\{i=1, \ldots, q\}.
    \]
\item[(A3)] The lower-level second order sufficient condition (LSOSC)  is satisfied at $(\bar x, \bar y, \bar u)$, where $\bar u$ is a lower-level Lagrange multiplier associated to $(\bar x, \bar y)$; i.e.,  
    \[
   d^\top \nabla^2_{yy}\ell(\bar x, \bar y, \bar u)d > 0 \quad \mbox{ for all } d\neq 0 \mbox{ s.t. }\; \left\{\begin{array}{rl}
       \nabla_y g_i(\bar x, \bar y)d=0 & i\in I(\bar x, \bar y), \\
       \nabla_y g_i(\bar x, \bar y)d=0  & i\in J(\bar u),
   \end{array}\right.
    \]
where $\ell$ is the lower-level Lagrangian function defined by 
\begin{equation}\label{eq:Lower-level-Lagrangian}
    \ell(x, y, u):= f(x, y) + u^\top g(x, y)
\end{equation}
 and the index set $J(\bar u)$ is given by  
\[
J(\bar u):=\left\{i\in[q]\;|\; \bar u_i >0\right\}.
\]
\item[(A4)] The lower-level strict complementary slackness (LSCS) is satisfied at the point $(\bar x, \bar y, \bar u)$; i.e., it holds that $I(\bar x, \bar y) = J(\bar u)$.
\end{itemize}

Under assumptions (A1)---(A4), the lower-level optimal solution function $y(\cdot)$ is well-defined and continuously differentiable from a neighborhood of $\bar x$ to a neighborhood of $\bar y$, and its Jacobian can be written as follows: 
\begin{equation}\label{eq:nablaY(x)}
   \begin{array}{rll}
    \nabla y(\bar x) &=& -\left(\nabla^2_{yy}\ell\right)^{-1}\Big\{\text{I}  \\[2ex]
                  & & -\left.\nabla_y g^\top_{\bar{I}}\left[\nabla_y  g_{\bar{I}}\left(\nabla^2_{yy}\ell\right)^{-1}\nabla_y g^\top_{\bar{I}}\right]^{-1} \nabla_y g_{\bar{I}}\left(\nabla^2_{yy}\ell\right)^{-1}\right\} \nabla^2_{xy}\ell\\[2ex]
    & &  -\left(\nabla^2_{yy}\ell\right)^{-1}\nabla_y g^\top_{\bar{I}}\left[\nabla_y g_{\bar{I}}\left(\nabla^2_{yy}\ell\right)^{-1}\nabla_y g^\top_{\bar{I}}\right]^{-1} \nabla_x g_{\bar{I}},
\end{array} 
\end{equation}
 where $\text{I}$ is the identity matrix of suitable dimensions and $\bar{I}\equiv I(\bar x, y(\bar x))$; see, e.g., Subsection 7.3.1 in the book \cite{shimizu2012nondifferentiable} (or \cite[Chapter 2] {fiacco1983introduction}). Furthermore, for full clarity, note that in the formula \eqref{eq:nablaY(x)}, we have 
$
 \ell\equiv \ell(\bar x, \bar y, \bar u) \mbox{ and } g_{\bar{I}}\equiv \left(g_i(\bar x, \bar y)\right)_{i\in \bar{I}}.
$

\begin{algorithm}
    \caption{Gradient descent algorithm with lower-level constraints}
    \label{alg:2}
    \begin{algorithmic}[1]
    \REQUIRE $x^0$ and $\{\alpha_k\}$;\\[2ex]
    \FOR{$k=0, 1, \dots, K$}
\STATE  $y(x^k)=\arg\underset{y}{\min}\left\{f(x^k, y)\left|\;\; y\in Y(x^k)\right.\right\}$;\\[2ex]
\STATE $I_k=I\left(x^k, y(x^k)\right)$;\\[2ex]
\STATE  Calculate $\nabla y(x^k)$ using \eqref{eq:nablaY(x)};\\[2ex]
\STATE  $\nabla\mathcal{F}(x^k)=\nabla_x F(x^k, y(x^k)) + \nabla y(x^k)^\top \nabla_y F(x^k, y(x^k))$;\\[2ex]
\STATE $x^{k+1}=x^k - \alpha_k\nabla\mathcal{F}(x^k)$\\[2ex]
\ENDFOR
    \end{algorithmic}
\end{algorithm}
 
 It goes without saying that under this framework, problem \eqref{eq:Implicit-Function-Model} is  well-defined and Algorithm \ref{alg:1} can be extended accordingly to this version of the problem with lower-level feasible set \eqref{eq:Y(x)_g(x,y)}, as it can be seen in Algorithm \ref{alg:2}.

Clearly, with the lower-level constraint, many more layers of difficulty appear in this extension of Algorithm \ref{alg:1}. First, at each iteration, the active indices of the current iteration point $(x^k, y(x^k))$ are needed; cf. line 3. Secondly, the complexity of calculating $\nabla y(x^k)$ increases significantly, as at each iteration. The formula requires computing the inverses of both $\nabla^2_{yy}\ell$ and $\nabla_y  g_{\bar{I}}\left(\nabla^2_{yy}\ell\right)^{-1}\nabla_y g^\top_{\bar{I}}$, as well as many other matrix operations, which need much more computing effort  in line 4  of Algorithm \ref{alg:2}. This is probably the key reason why almost all bilevel learning algorithms avoid lower-level constraints. It might be useful to observe that if $g\equiv 0$ in Algorithm \ref{alg:2}, the method just reduces to Algorithm \ref{alg:1}. 

Besides the challenge in accommodating constraints in the lower-level problem, we would like to draw attention to the fact that even when $y(\cdot)$ is well-defined as a vector-valued function in some neighborhood of a point of interest, it is nonsmooth in general. The premise of having $y(\cdot)$ as a continuously differentiable function in some neighborhood of interest requires all the functions involved in the lower-level problem to be at least twice continuously differentiable, as required in almost all bilevel learning papers (also recalled this in the context of constraints above). However, the loss functions in a wide range of machine learning tasks (and by extension the lower-level objective functions of multiple bilevel learning applications) are nonsmooth; see, e.g., \cite{terven2025comprehensive,wang2022comprehensive} or \cite{liu2023rethinking,shi2022efficient,okuno2021,alcantara2021unified}. To restore smoothness of lower-level functions, some bilevel learning papers have applied smoothing techniques (see, e.g., \cite{alcantara2021unified,okuno2021}) or use some usual tricks that result in possibly introducing constraints in problems that are initially unconstrained (see, e.g., \cite{bennett2006model,Kunapuli2008ClassificationMS,ward2025mathematical}).

\begin{figure}[H]
    \centering
    \includegraphics[scale = 0.45]{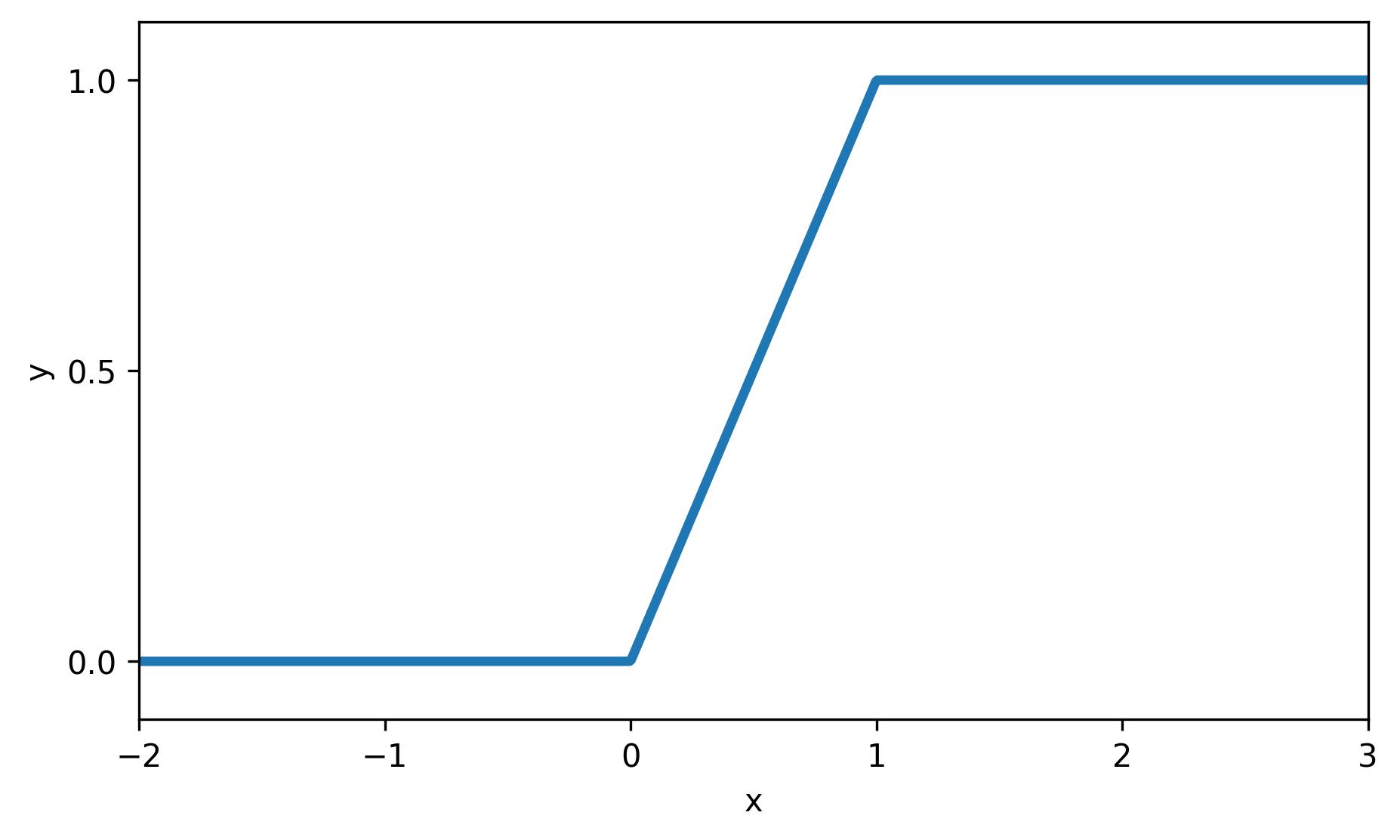}
    \caption{Graph of the optimal solution function $y(\cdot)$ for the example in \eqref{example-illus-1}.}
    \label{fig:y() example}
\end{figure}

Moreover, it is important to observe that having twice continuously differentiable functions in the lower-level problem does not necessarily guarantee that $y(\cdot)$ is smooth. For instance, for the toy lower-level problem example 
\begin{equation}\label{example-illus-1}
    \min_y \left\{\left.\frac{1}{2}(y-x)^2\,\right|\;\;  0\leq y\leq 1\right\},
\end{equation}
the LLICQ (i.e., (A2)) and LSOSC (i.e., (A3)) are satisfied at $(0, 0)$ and $(0, 0, 0)$, respectively, as $\bar u =0$ here. The graph of $y(\cdot)$ for this example can be seen in Figure \ref{fig:y() example}. There, you can see that $y(\cdot)$ is nondifferentiable at $0$ because the LSCS fails, given that 
\[
I(0, 0)=\{1\} \; \mbox{ and }\; J(0)=\emptyset \;\,(\mbox{with }\, g_1(x,y):=-y \, \mbox{ and } \, g_2(x,y):=y-1).
\]
As it can be seen in Figure \ref{fig:y() example}, $y(\cdot)$ is a piecewise smooth function in the context of this example. This is a rather general observation for a constrained parametric optimization problem. 
To be more precise, consider the lower-level problem present in  \eqref{eq:Bilevel-Optimization-Problem}, where the lower-level feasible set is described by \eqref{eq:Y(x)_g(x,y)}. If assumptions (A1)---(A4) are satisfied at some point $(\bar x, \bar y)$,
then, there are open neighborhoods $U$ of $\bar x$ and $V$ of $\bar y$, and a continuous function $y(\cdot)$ mapping $U$ to $V$ such that for each $x\in U$, $y(x)$ is the unique local solution of the lower-level problem (for $x$ fixed) in $V$; and the function $y(\cdot )$ is piecewise continuously differentiable around $\bar x$. This means that $y(\cdot)$ is continuous and there is a finite family of continuously differentiable functions $y^1(x)$, \ldots, $y^N(x)$ defined in a neighborhood of $\bar x$, such that $y(x)\in \{y^l(x), \ldots, y^N(x)\}$ for any $x\in U$; see \cite{ralph1995directional,dempe2002foundations}. Under this framework, $y(\cdot)$ is local Lipschitz continuous near $\bar x$, as it is the case for any piecewise smooth function \cite{scholtes2012introduction}.

So, under relatively affordable assumptions, $y(\cdot)$ is a locally Lipschitz continuous function. There are multiple works on approaches to calculate generalized derivatives of $y(\cdot)$ (see \cite{fiacco1983introduction,ralph1995directional,dempe2002foundations} and references therein) and their use in the general field of bilevel optimization; see, e.g., \cite{kolstad1990derivative,savard1994steepest,vicente1994descent}.  However, these results do not seem to have been implemented in the context bilevel learning; and one of the main reasons for this is that these generalized derivatives of $y(\cdot)$ are generally abstract in nature or very difficult to calculate in practice. Recently, the paper \cite{bolte2021conservative} introduced the concept of \textit{conservative derivative} to efficiently extend derivative approximation tools (such as \textit{automatic differentiation}, a backbone object to enhance deep learning algorithms \cite{baydin2018automatic}) to nonsmooth functions. 
The concept of conservative derivative has been used in \cite{grazzi2024nonsmooth} to implement the classical iterative differentiation (ITD) and approximate implicit differentiation (AID) schemes, widespread in the practical implementation of lines 3 and 4 of Algorithm \ref{alg:1} in the smooth unconstrained-lower-level bilevel learning problem, to the case where $y(\cdot)$ is a nonsmooth function. More work is needed to generalized derivative approximation schemes to broader classes of BL problems. 


\subsection{Restoring the implicit function approach for problems with multiple lower-level optimal solutions}
The biggest challenge for the state of the art method for bilevel learning, which is sketched in Algorithm \ref{alg:1}, is  the requirement to always have a unique solution for the lower--level problem; cf. \eqref{eq:S(x)=1}. If this assumption fails, the implicit function approach \eqref{eq:Implicit-Function-Model} is out of question. Ensuring that this condition holds requires very strong assumptions as we have discussed above. Unfortunately, this assumption cannot hold for a wide range of BL applications; see, e.g., \cite{liu2022general,liu2020generic,sow2022primal,yao2024constrained,xiao2023generalized,guan2022gradient,huo2024perturbed}, for a selection of such problems that do not satisfy condition \eqref{eq:S(x)=1}. To address this challenge, a series of articles from the BO literature have proposed a regularization approach, which consists to introduce a perturbation to the lower-level problem of the form
 \begin{equation}\label{eq:S_alpha(x)}
    \underset{y\in Y(x)}{\text{min}}~f(x,y) + \alpha \psi(x, y),
\end{equation}
in order to force the fulfillment of condition \eqref{eq:S(x)=1}. In fact, if the lower-level problem is \textit{just} convex and the function $\psi : \mathbb{R}^n\times \mathbb{R}^m \rightarrow \mathbb{R}$ is strongly convex in $y$, then it holds that 
\[
\left\{\left. x\in X\right|\;\; |S_{\alpha}(x)| \neq 1 \right\} = \emptyset \;\mbox{ for all }\; \alpha >0,
\]
where $S_{\alpha}$ is the optimal solution mapping for the regularized problem \eqref{eq:S_alpha(x)}; i.e., 
\[
S_{\alpha}(x) := \underset{y\in Y(x)}{\text{argmin}}~f(x,y) + \alpha \psi(x, y).
\]
In this context, model \eqref{eq:Implicit-Function-Model} can be rescued as follows
\begin{equation}\label{eq:Implicit-Function-Model-alpha}
      \underset{x\in X}{\min}~\mathcal{F}_{\alpha}(x):=F(x, y_{\alpha}(x)) \;\;\mbox{ with }\;\; \{y_{\alpha}(x)\}=S_{\alpha}(x) \;\mbox{ for all }\; \alpha>0.
\end{equation}
With the additional caution ensuring that $y_{\alpha}(\cdot)$ describes a smooth function for all $\alpha >0$, Algorithm \ref{alg:1} can be extended to problems with multiple optimal solutions if $Y(x)=\mathbb{R}^m$, and similarly, for problems with lower-level constraints described by \eqref{eq:Y(x)_g(x,y)}, a suitable version of Algorithm \ref{alg:2} can also be extended to this case. 

To make such an iterative procedure rigorous, a sequence $(\alpha_n)_n$ such that $\alpha_n >0$ for all $n\in \mathbb{N}$ with $\alpha_n \downarrow 0$ can be introduced with 
\begin{equation}\label{eq:Implicit-Model-Sequence-Inclusion}
x(\alpha_n) \in \underset{x\in X}{\text{argmin}}~\mathcal{F}_{\alpha_n}(x) \;\mbox{ for }\; n\in \mathbb{N}.  
\end{equation}
The regularization in \eqref{eq:S_alpha(x)}, due to Tikhonov \cite{tikhonov1965regularisation}, has been used in \cite{dempe1996algorithm} to develop a gradient descent-type algorithm, with $\psi(x, y):= F(x, y)$ if $\nabla^2_{yy} F(x, y)$ is positive definite for all $(x, y)$, to solve a constrained lower-level bilevel program. On the other hand, the regularization term  $\psi(x, y):= \|y\|^2$ is used in \cite{dempe2000bundle,dempe2001bundle} to develop bundle algorithms for different classes of the bilevel optimization with lower-level constraints; further details on regularization methods for bilevel programs can be found in \cite{dempe2002foundations}. 


Despite suitable technical assumptions in the aforementioned papers to ensure the convergence of such algorithms, very simple examples can be constructed to show that the resulting limit point $\bar x = \lim_{n\rightarrow \infty} x(\alpha_n)$ can be quite far away from the true optimal solution or stationary point of the corresponding version of problem \eqref{eq:Bilevel-Optimization-Problem}; see, e.g., \cite{dempe1996algorithm,morgan2006stackelberg}. However, it can be shown that under mild assumptions, the limit point $(\bar x, y(\bar x))$ resulting from an algorithm based on the Tikhonov regularization above converges to a \textit{lower Stackelberg equilibrium}; i.e., a point that satisfies
\[
\bar x\in X, \;\, \bar y\in S(\bar x), \quad F(\bar x, \bar y) \leq \underset{x\in X}{\inf}\,\underset{y\in S(\bar x)}{\sup} F(x, y).
\]
This obviously implies that a lower Stackelberg equilibrium lies between optimistic and pessimistic optimal solution of problem \eqref{eq:Bilevel-Optimization-Problem} given that 
\[
\underset{x\in X}{\inf}\,\underset{y\in S(\bar x)}{\inf} F(x, y) \leq \;\, F(\bar x, \bar y)\;\,  \leq \underset{x\in X}{\inf}\,\underset{y\in S(\bar x)}{\sup} F(x, y).
\]
What this could mean in bilevel learning is that the Tikhonov regularization approach above could be a tractable framework to extend Algorithms \ref{alg:1} and \ref{alg:2} to problems with multiple lower-level optimal solutions if one is inclined to accept lower Stackelberg equilibrium. Also note that the concept of lower Stackelberg equilibrium is closely connected to the notion of \textit{subgame perfect Nash equilibrium}, introduced by the Nobel laureate Reinhard Selten in \cite{selten1965spieltheoretische} and which is widely used in economics. 
%

\subsection{Implicit set-valued model}
A direct approach to deal with the failure of condition \eqref{eq:S(x)=1} could simply be to insert the lower-level optimal solution mapping $S$ in the upper-level objective function, therefore leading to the set-valued optimization problem 
\begin{equation}\label{eq:Set-valued-Implicit-Model}\tag{\mbox{$\text{P}_{\text{S}}$}}
      \underset{x\in X}{\min}~\mathcal{F}_S(x):=F(x, S(x)).
\end{equation}
This model, which is a direct extension of the implicit function problem \eqref{eq:Implicit-Function-Model} to bilevel optimization problem where the  lower-level optimal solution mapping is set-valued, was first studied in \cite{zemkoho2016solving} using the corresponding Pareto optimal solution concept. Namely, in the latter paper, optimality conditions for problem \eqref{eq:Set-valued-Implicit-Model} were derived and shown to capture all the stationary notions known for the optimistic problem \eqref{eq:Optimistic-Model}. More recently, it was shown in \cite{som2025bilevel} that problem \eqref{eq:Set-valued-Implicit-Model} can be equivalent to the optimistic problem \eqref{eq:Optimistic-Model} (resp. pessimistic problem \eqref{eq:Pessimistic-Model}), while considering the l-minimal optimal solution (resp. u-minimal optimal solution) for problem \eqref{eq:Set-valued-Implicit-Model}, in the sense of set-valued optimization, under mild assumptions. 

It must be said that so far, the model \eqref{eq:Set-valued-Implicit-Model} is only a theoretical object, which can  enable some enhanced insights on the bilevel optimization problem when the lower-level player has multiple options to make their decision for some choices of the upper-level player. However, there is not much progress in solving set-valued optimization problems in the current literature. For some recent attempts to compute optimal solutions for the problem, interested readers are referred to \cite{bouza2021steepest,ghosh2025quasi,lohne2025solution}, for example.


\section{Karush-Kuhn-Tucker reformulation-based methods}\label{sec:Constrained optimization}
We have seen in the previous section that when the lower-level problem is constrained, a direct extension of the classical bilevel learning algorithm presented in Section \ref{sec:Main algorithmic techniques} is intractable; cf. Algorithm \ref{alg:2}. Throughout this section, we continue with the assumption that the lower-level problem is constrained by the set \eqref{eq:Y(x)_g(x,y)}. However, unlike in the previous section, where we need the very strong assumptions (A1)--(A4) for the corresponding version of \eqref{eq:Implicit-Function-Model} to be well-defined as a smooth optimization problem, here, we only need to impose that the lower-level problem is convex (i.e., for all $x\in X$, the functions $f(x,.)$ and $g_i(x, .)$, for $i=1, \ldots, q$, are convex) and satisfies a constraint qualification (CQ) that can enable us to write the Karush-Kuhn-Tucker (KKT) conditions that characterize the inclusion $y\in S(x)$. 

As we will refer to the Mangasarian-Fromovitz constraint
qualification (MFCQ) multiple times in the sequel, note that for a general constrained optimization problem 
\begin{equation}\label{eq:min_standard}
    \begin{array}{rl}
    \min &\mathfrak{f}(x)\\[1ex]
    \mbox{s.t.} & x\in \mathcal{C}:=\left\{\left. x\in \mathbb{R}^n\right|~\mathfrak{g}(x)\leq 0, \;\; \mathfrak{h}(x)=0\right\}
    \end{array}
\end{equation}
with continuously differentiable functions $\mathfrak{f} : \mathbb{R}^n \rightarrow \mathbb{R}$, $\mathfrak{g} : \mathbb{R}^n \rightarrow \mathbb{R}^p$, and $\mathfrak{h} : \mathbb{R}^n \rightarrow \mathbb{R}^q$, it will be said to hold at a point $\bar x\in \mathcal{C}$ if the following condition is satisfied:
\begin{equation}\label{eq:MFCQ}\tag{MFCQ}
    \left.\begin{array}{r}
        \nabla\mathfrak{g}(\bar x)^\top \alpha + \nabla\mathfrak{h}(\bar x)^\top \beta =0\\[1ex]
         \alpha\geq 0, \;\; \alpha^\top\mathfrak{g}(\bar x)=0
    \end{array}\right\} \Longrightarrow \left\{\begin{array}{l}
         \alpha =0,\\[1ex]
         \beta =0.
    \end{array} \right.
\end{equation}

Throughout this section, we assume here that the lower-level problem is convex and \eqref{eq:MFCQ} holds at all $(x, y)\in \text{gph}\,Y$ with $x\in X$ and $y\in S(x)$, then, as mentioned in Section \ref{sec:A short history}, historically, the basic approach to solve the corresponding standard optimistic bilevel program \eqref{eq:Standard-Optimistic-Bilevel} has consisted to first write it as a single-level optimization problem  
\begin{equation}\label{eq:KKT-Optimistic-Bilevel}\tag{KKT}
\begin{array}{rl}
\underset{x, y, u}{\min} & F(x, y)\\
\mbox{s.t.} & x\in X, \;\,  \nabla_y \ell(x, y, u)=0,\\[1ex]
            & u\geq0, \; g(x,y)\leq 0, \; u^\top g(x,y)=0,
\end{array}
\end{equation}
known as KKT reformulation. Here, the lower-level Lagrangian function $\ell$ is defined as in \eqref{eq:Lower-level-Lagrangian}. 
Problem \eqref{eq:KKT-Optimistic-Bilevel} has been one of the main \textit{go to} frameworks to develop numerical methods for the optimistic bilevel optimization problem \eqref{eq:Standard-Optimistic-Bilevel} since its introduction in the field of mathematical optimization; see, e.g., \cite{dempe2002foundations,Colson2007AnOO,dempe2020bilevel,dempe2012karush} and references therein. 

After some discussion on the background the KKT reformulation problem \eqref{eq:KKT-Optimistic-Bilevel}, we provide an overview of classical ideas from the BO literature that remain mostly unexplored in BL. At the end of the section we present some works on BL that have applied the KKT reformulation. 
Overall, the material here could serve as base for investigations on the scalability of these techniques, especially in terms of their approximations with ideas similar to the ones that have been at the core of the success of the gradient descent method for bilevel learning described in Section \ref{sec:Main algorithmic techniques}. 

The basic idea behind reformulation \eqref{eq:KKT-Optimistic-Bilevel} is the fact that under the lower-level convexity assumption and the fulfillment of the LLICQ at a point $(x, y)\in \text{gph}\,S$, 
\begin{equation}\label{eq:LL_KKT_System}
    y\in S(x) \Longleftrightarrow \exists u\in \mathbb{R}^q: \left\{\begin{array}{l}
   \nabla_y \ell(x, y, u)=0,\\
    u\geq0, \; g(x,y)\leq 0, \; u^\top g(x,y)=0.
\end{array}\right.
\end{equation}
Despite this equivalence, the relationship between problems \eqref{eq:Standard-Optimistic-Bilevel} and \eqref{eq:KKT-Optimistic-Bilevel} is a bit tricky due to the appearance of the Lagrange multiplier  $u$ in the new problem. Both problems are globally equivalent in some sense \cite{dempe2012bilevel}. However, locally, to get a local optimal solution $(\bar x, \bar y)$ of problem \eqref{eq:Standard-Optimistic-Bilevel} from \eqref{eq:KKT-Optimistic-Bilevel}, one needs to ensure that $(\bar x, \bar y, u)$ is locally optimal for the latter problem for all $u\in \Lambda(\bar x, \bar y)$, where
\[
\Lambda(\bar x, \bar y):=\left\{\left. u\in \mathbb{R}^q\,\right|\;  \nabla_y \ell(\bar x, \bar y, u)=0, \; u\geq0, \; g(\bar x, \bar y)\leq 0, \; u^\top g(\bar x, \bar y)=0\right\}
\]
denotes the set of lower-level Lagrange multipliers corresponding to the lower-level optimal solution $\bar y$ when $x$ is fixed at $\bar x$. Given that the set $\Lambda(\bar x, \bar y)$ can be made of infinitely many points, this is an intractable prospect, even for very small toy examples, not to imagine this in the context of BL. For more details on this connection between problems \eqref{eq:Standard-Optimistic-Bilevel} and \eqref{eq:KKT-Optimistic-Bilevel}, see \cite{dempe2012bilevel}. For a systematic analysis of single-level reformulations of problems \eqref{eq:KKT-Optimistic-Bilevel} based on transformations of the lower-level problem (especially from the duality theory perspective), interested readers are referred to  \cite{dempe2025duality}. 

As a transition point to discuss methods to solve problem \eqref{eq:KKT-Optimistic-Bilevel}, it would be useful to note that some care is needed in handling the problem, and not just treat it as any random mathematical program with equality and inequality constraints. For example, if we do so, the first issue that we are faced with is that many classical constraint qualifications, including the MFCQ, will fail; see, e.g.,  \cite{scheel2000mathematical}.  
The challenges in solving problem \eqref{eq:KKT-Optimistic-Bilevel} are due to the complementarity conditions, present in the last line of the feasible set, and which do not allow for a feasible point to strictly satisfy the inequality constraints to exist. This problematic structure of the feasible set of problem \eqref{eq:KKT-Optimistic-Bilevel} has motivated the development of specially tailored algorithms to solve it. 
Throughout the literature, there are roughly four algorithmic techniques to tackle the problem, and we briefly outline them below.  

\subsection{Partial penalization--based methods} Partial penalization consists of penalizing the term $u^\top g(x,y)$ by moving it from the feasible set to the (upper-level) objective function. This results in the new problem 
\begin{equation}\label{eq:lambda-KKT-Optimistic-Bilevel}\tag{\mbox{KKT$_\lambda$}}
\begin{array}{rl}
 \underset{x, y, u}{\min} &F(x, y)-\lambda u^\top g(x,y) \\
 \mbox{s.t.} & x\in X, \;\,  \nabla_y \ell(x, y, u)=0, \;\, u\geq0, \; g(x,y)\leq 0
\end{array}
\end{equation}
with penalization parameter $\lambda >0$. The feasible set becomes much easier to handle with transformation \eqref{eq:lambda-KKT-Optimistic-Bilevel}, as the \eqref{eq:MFCQ}, for example, can usually be satisfied for this penalized  problem. 
Problem \eqref{eq:lambda-KKT-Optimistic-Bilevel} has been widely used since the early days, as a based to solve problem \eqref{eq:KKT-Optimistic-Bilevel}, and therefore the corresponding bilevel optimization \eqref{eq:Standard-Optimistic-Bilevel}. It was possibly first used to  tackle the fully linear bilevel program in \cite{campelo2000note} (with improvements later provided in \cite{white1993penalty}), where the penalty parameter is sequentially increased until a stopping criterion is achieved. 
The penalization model \eqref{eq:lambda-KKT-Optimistic-Bilevel} with a sequentially varying penalization parameter was also used in \cite{leyffer2006interior} to develop an interior-point-type method for the general MPEC problem. For such a class of problem, the partial penalization model of the form \eqref{eq:lambda-KKT-Optimistic-Bilevel} is thoroughly studied in \cite{ralph2004some} but in the case of a fixed parameter. 

Partial penalization has also been used in \cite{ye1997exact} (also see \cite{ralph2004some}) as a form of constraint qualification to derive  necessary optimality conditions for problem \eqref{eq:KKT-Optimistic-Bilevel}. More recently, detailed algorithmic studies on the partial penalization approach in \eqref{eq:lambda-KKT-Optimistic-Bilevel} have been conducted in \cite{zemkoho2021theoretical,ward2025mathematical} for the standard optimistic bilevel optimization problem \eqref{eq:Standard-Optimistic-Bilevel}. As it is common in exact penalization methods, finding a good value for the parameter $\lambda$ for problem \eqref{eq:lambda-KKT-Optimistic-Bilevel} is difficult. The paper \cite{zemkoho2021theoretical} also includes an empirical study on the subject. 

\subsection{Relaxation--based methods} They consist of a process that starts with the enlargement of the feasible set of problem \eqref{eq:KKT-Optimistic-Bilevel} by introducing a relaxation function to generate a more tractable subproblem. The standard relaxation schemes for problem   \eqref{eq:KKT-Optimistic-Bilevel} can be summarized in the compact model
\begin{equation}\label{eq:t-KKT-Optimistic-Bilevel}\tag{\mbox{\text{KKT$_t$}}}
\begin{array}{rl}
\underset{x, y, u}{\min} & F(x, y)\\
\mbox{s.t.} & x\in X, \;\,  \nabla_y \ell(x, y, u)=0,\;\; \phi^t_{i, \mathcal{R}}(x, y, u)\leq 0, \; i\in[q].
\end{array}
\end{equation}
Here, $t>0$ denotes the relaxation parameter, while $\mathcal{R}$ is used to represent a specific relaxation from the literature; more precisely,   $\mathcal{R}\in \left\{\mbox{S}, \mbox{LF}, \mbox{KDB}, \mbox{SU}, \mbox{KS}\right\}$, where S, LF, KDB, SU, and KS respectively represents the Scholtes, Lin and Fukushima, Kadrani, Dussault and Benchakroun, Steffensen and Ulbrich, and Kanzow and Schwartz relaxation of problem \eqref{eq:KKT-Optimistic-Bilevel}. 
For $t>0$ and $i=1, \ldots, q$, the function $\phi^t_{i, \mathcal{R}}$ can be defined by  
\begin{equation}
    \phi^t_{i, \mathcal{R}}(x, y, u):=\left\{ 
    \begin{array}{lll}
       \left(\begin{array}{cc}
            g_i(x,y)  \\
            -u_i\\
            -u_i g_i(x,y)-t
       \end{array} \right)  &  \mbox{ if } & \mathcal{R}:=S,\\[1.5ex]
\left(\begin{array}{cc}
           -(u_{i}g_{i}(x,y)+t^{2}) \\[1ex]
-(u_{i}+t)(-g_{i}(x,y)+t)+t^{2}
       \end{array} \right)  &  \mbox{ if } & \mathcal{R}:=LF,\\[3ex]
\left(\begin{array}{cc}
            g_i(x,y) -t  \\
            -u_i - t\\
            -(u_i -t)(g_i(x,y) + t)
       \end{array} \right)  &  \mbox{ if } & \mathcal{R}:=KDB,\\[3ex]
 \left(\begin{array}{cc}
            g_i(x,y)  \\ 
            -u_i\\
            \varphi^t_{i, SU}(x,y,u)
       \end{array} \right)  &  \mbox{ if } & \mathcal{R}:= SU,\\[3ex]
 \left(\begin{array}{cc}
            g_i(x,y)  \\
            -u_i\\
             \varphi^t_{i, KS}(x,y,u)
       \end{array} \right)  &  \mbox{ if } & \mathcal{R}:=KS,
    \end{array}
    \right.
\end{equation}
where, for $t>0$ and $i=1, \ldots, q$, the function $\varphi^t_{i, SU}$ is defined as 
{\small{
\[
\varphi^t_{i, SU}(x,y,u):=\left\{
\begin{tabular}{lll}
$2u_{i}$ & if & $g_{i}(x,y)+u_{i}\leq-t$,\\
$-2g_{i}(x,y)$ & if & $g_{i}(x,y)+u_{i}\geq t$,\\
$u_{i}-g_{i}(x,y)-t\theta\left(\dfrac{u_{i}+g_{i}(x,y)}{t}\right)$ & if & $\left\vert
u_{i}+g_{i}(x,y)\right\vert <t.$
\end{tabular}
\right.
\]
}
}
As for the function $\theta(\cdot)$, it represents a suitable regularization function (see details in \cite{steffensen2010new}) 
and $\varphi^t_{i, KS}$ is defined by
\[
\varphi^t_{i, KS}(x,y,u):=\left\{
\begin{tabular}{ll}
$(u_{i}-t)(-g_{i}(x,y)-t)$ & $\text{if\;\; } u_{i}-g_{i}(x,y)\geq2t$,\\[1ex]
$-\dfrac{1}{2}\left((u_{i}-t)^{2}+(-g_{i}(x,y)-t)^{2}\right)$ & $\text{if\;\; }u_{i}
-g_{i}(x,y)<2t.$
\end{tabular}
\right.
\]
For an overview and detailed study of these relaxations in the broader context of mathematical programs with complementarity constraints, interested readers are referred to \cite{hoheisel2013theoretical} and references therein, where specific advantages and drawbacks of each relaxation are given and compared. It might be useful to highlight that the S-relaxation, which is probably the first and simplest one, introduced in \cite{scholtes2001convergence}, seems to usually show the best numerical performance (even tough convergence is typically to a C-stationary point, which is relatively weak). 

Note that the complementarity conditions present in \eqref{eq:KKT-Optimistic-Bilevel} make the feasible set of the problem very thin (as it is the union of segments of the axes), and hence, the process for a numerical procedure to search for an optimal point from it to be quite tricky. For a given relaxation above, the parameter $t>0$ helps to control the enlargement of the feasible set such that as $t\downarrow 0$, one generally recaptures the feasible set of problem \eqref{eq:KKT-Optimistic-Bilevel}.  These relaxation methods typically lead to C or M-stationarity points, depending on the assumptions made to establish their convergence results. For more details on the related theory and the definitions of these stationarity concepts,  interested readers are referred to the article \cite{hoheisel2013theoretical}; if one is interested in the implementation of these relaxations in the context of the pessimistic bilevel program \eqref{eq:Pessimistic-Model}, see \cite{benchouk2025scholtes,benchouk2024relaxation}, for example.

\subsection{NCP function--based methods} These are algorithms to solve problem \eqref{eq:KKT-Optimistic-Bilevel} that proceed by first transforming the complementarity conditions into a system of equations in such a way to get an equivalent problem of the form  
\begin{equation}\label{eq:NCP-Optimistic-Bilevel}\tag{NCP}
\begin{array}{rl}
\underset{x, y, u}{\min} & F(x, y)\\
\mbox{s.t.} & x\in X, \;\,  \nabla_y \ell(x, y, u)=0,\;\; \phi_{i, \mathcal{N}}(x, y, u)= 0, \, i\in[q],
\end{array}
\end{equation}
where $\phi_{i, \mathcal{N}}$, $i=1, \ldots, q$,  represents a so-called nonlinear complementarity problem
(NCP) function. Precisely, for a given index $i=1, \ldots, q$, an NCP function $\phi_{i, \mathcal{N}}$ is a real-valued function constructed such that we have 
\begin{equation}\label{NCP-function-i}
    \phi_{i, \mathcal{N}}(x, y, u)= 0 \Longleftrightarrow \;\, \left[u_i\geq0, \; g_i(x,y)\leq 0, \; u_i g_i(x,y)=0 \right].
\end{equation}
The concept of NCP function, introduced by Mangasarian in \cite{mangasarian1976equivalence}, has been widely studied in the literature, considering the occurrence of complementarity conditions in many practical applications in areas such as economics, engineering, and science, just to name a few. The most famous NCP functions are probably the \textit{Fischer--Burmeister} and \textit{min} operator functions, which are respectively defined from $\mathbb{R}^2$ to $\mathbb{R}$ by
\[
\phi_{\text{FB}}:=\sqrt{a^2 + b^2}-(a+b) \;\mbox{ and }\; \phi_{\min}(a, b):=\min\{a, \, b\},
\]
which each vanishes at a point $(a, b)\in \mathbb{R}^2$ if and only if $a\geq 0$, $b\geq 0$, and $ab=0$. Note that using this NCP function in \eqref{NCP-function-i}, we have $\phi_{i, \mathcal{N}}(x, y, u):=\phi_{\text{FB}}\left(-g_i(x,y),\; u_i\right)$ for $i\in [q]$. The Fischer-Burmeister function, introduced in \cite{fischer1992special}, has been particularly prominent due to the fact that the associated merit function 
$\psi(a, b):=\frac{1}{2}\|\phi_{\text{FB}}(a, b)\|^2$ for the system $\phi_{\text{FB}}(a, b)=0$ is continuously differentiable. This has enable the development of powerful algorithms for semismooth systems of equations involving complementarity conditions; see, e.g.,  \cite{de1996semismooth} for an important optimization algorithm in this context. 

Note that there are so many ways to construct NCP functions; see, e.g., \cite{alcantara2020construction,galantai2012properties} for some recent studies on the subject. With regards to solving  \eqref{eq:KKT-Optimistic-Bilevel}, the NLPEC solver  (\href{https://www.gams.com/50/docs/S\_NLPEC.html}{https://www.gams.com/50/docs/S$\_$NLPEC.html}) works by automatically enabling the selection of an NCP function for a given general MPEC problem, to design a process to compute solutions for problem \eqref{eq:NCP-Optimistic-Bilevel}. Recently, special NCP functions that improve the development of efficient methods for neural network computations have been discovered; see, e.g., \cite{alcantara2019neural,xie2021neural}.
It must however be said that a detailed study of the impact of NCP functions in the numerical development of methods for the KKT reformulation \eqref{eq:KKT-Optimistic-Bilevel} for the optimistic bilevel optimization problem \eqref{eq:Standard-Optimistic-Bilevel} has not yet been done. 

\subsection{Nonsmooth equation system--based methods} 
Given that  \eqref{eq:KKT-Optimistic-Bilevel} is a nonconvex optimization problem all the methods discussed so far to solve it can typically be shown to theoretically only converge to stationary points. However, considering the nature of the feasible set of the problem, defined by complementarity conditions, there are multiple different types of stationarity concepts, depending on the specific approach used or the properties imposed for the problem data to ensure theoretical convergence. The main types of stationarity concepts for problem \eqref{eq:KKT-Optimistic-Bilevel} are the C--, M--, and S--stationarity concepts, where the latter is the \textit{strong} stationarity concept, which is equivalent to the KKT conditions of problem \eqref{eq:KKT-Optimistic-Bilevel} when it is viewed as a usual optimization problem with equality and inequality constraints. As for M--and C--, they  stand for the \textit{Mordukhovich} and \textit{Clarke} stationarity conditions, respectively, due to the variational analysis tools applied to compute the generalized derivative of the involved nonsmooth functions; for more details on these concepts and how to derive them for local optimal solutions of \eqref{eq:KKT-Optimistic-Bilevel}, see, e.g., \cite{scheel2000mathematical,flegel2005m,dempe2012karush}.

Another important common point of the methods described so far for problem \eqref{eq:KKT-Optimistic-Bilevel} is that they all rely on first building a \textit{nice} auxiliary problem that is subsequently solved with the hope compute a \textit{solution} of the original problem. As we have just said above, theretically, these algorithms can only be shown to compute stationary points. However, another classical philosophy in solving constrained optimization problems, with continuous variables, is to directly compute these stationary points. This has led to the stream of work on semismooth Newton methods around the early 1990s for general smooth constrained problems; see, e.g., \cite{fischer1992special,qi1993nonsmooth,pang1993nonsmooth}. This class of methods has recently been explored to compute M-stationarity points for mathematical programs with complementarity constraints, and can therefore be used to tackle the KKT reformulation \eqref{eq:KKT-Optimistic-Bilevel}; see, e.g., \cite{guo2015solving,harder2021reformulation,wu2015inexact}. 

To get a flavor of how a nonsmooth equation-based method can work in practice, consider the KKT reformulation based Lagrangian function $L_{\text{K}}$ defined for any $(x, y, u, \alpha, \beta, \gamma)$ with $(x, y, u)\in \mathbb{R}^n\times \mathbb{R}^m \times \mathbb{R}^q$ and $(\alpha, \beta, \gamma)\in \mathbb{R}^p\times \mathbb{R}^q \times \mathbb{R}^m$ by
\[
L_{\text{K}}(x, y, u, \alpha, \beta, \gamma):= F(x, y) + \alpha^\top G(x) + \beta^\top g(x,y) + \gamma^\top \ell(x, y, u),
\]
where, $G : \mathbb{R}^n \rightarrow \mathbb{R}^p$ is a continuously differentiable function that describes the upper-level feasible set as follows:
\begin{equation}\label{eq:X-G(x)}
    X:=\left\{x\in \mathbb{R}^n\left|\; G(x) \leq 0  \right.\right\}.
\end{equation}
Additionally, for any feasible point $(x, y, u)\in \mathbb{R}^n\times \mathbb{R}^m \times \mathbb{R}^q$ of problem \eqref{eq:KKT-Optimistic-Bilevel}, the classical partition of the index set associated to the complementarity conditions that partly describe the feasible set of the problem is given by 
\[
\begin{array}{l}
    \eta \, := \, \eta(x, y,  u) \, := \, \left\{i=1, \ldots, q\left|\;\,{u}_i=0,\;\; g_i(x,  y)>0 \right.\right\},\\[0.5ex]
     \mu \, := \, \mu(x, y,  u) \, := \, \left\{i=1, \ldots, q\left|\;\,{u}_i=0,\;\; g_i(x,  y)=0 \right.\right\},\\[0.5ex]
      \nu \, := \, \nu(x, y,  u) \, := \, \left\{i=1, \ldots, q\left|\;\,{u}_i>0,\;\; g_i(x, y)=0 \right.\right\}.
\end{array}
\]
Based on this notation, a feasible point $(x, y, u)\in \mathbb{R}^n\times \mathbb{R}^m \times \mathbb{R}^q$ of problem \eqref{eq:KKT-Optimistic-Bilevel} is said to be a M-stationary point if there exists Lagrange multipliers $\alpha\in \mathbb{R}^p$, $\beta\in \mathbb{R}^q$ and $\gamma\in \mathbb{R}^m$ such that the following necessary optimality conditions hold: 
\begin{eqnarray}
    \nabla_{x, y}L_{\text{K}}(x, y, \alpha, \beta, \gamma)=0, \label{M-Stat1}\\[0.5ex]
    \alpha \geq 0,\;\; G(x)\leq 0, \;\; \alpha^\top G(x)=0,\label{M-Stat2}\\[0.5ex]
    \nabla_y g_\nu(x,y)\gamma =0, \;\; \beta_\eta =0,\label{M-Stat3}\\[0.5ex]
    \forall i\in \mu: \;\; \left(\beta_i>0\,\wedge \, \nabla_y g_i(x,y)\gamma>0\right) \,\vee\, \left(\beta_i \nabla_y g_i(x,y)\gamma\right) =0. \label{M-Stat4}
\end{eqnarray}
For details on how to obtain these conditions for a given local optimal solution of problem \cite{dempe2012karush,dempe2013bilevel}.
The optimality conditions \eqref{M-Stat1}--\eqref{M-Stat4} can be written as a nonlinear system of equations if we consider the function the Fischer-Burmeister function (introduced in \cite{fischer1992special})
\begin{equation}\label{eq:phi_FB}
    \phi_{\text{FB}}(a, b):=\sqrt{a^2 + b^2}-(a+b) \;\mbox{ for }\; (a, b)\in \mathbb{R}^2
\end{equation}
and the M-stationarity function (introduced in \cite{harder2021reformulation}) 
\[
\phi_{\text{M}}(a, b, c, d):= \min\left\{\begin{array}{c}
     \max\left\{-a,\, |b|,\, |c|\right\},\\
     \max\left\{-b,\, |a|,\, |d|\right\},\\
     \max\left\{|a|,\, |b|,\, c,\, d\right\}
\end{array}\right\}
\;\mbox{ for }\; (a, b, c, d)\in \mathbb{R}^4.
\]
The M-stationarity system \eqref{M-Stat1}--\eqref{M-Stat4} can be equivalently written as 
\begin{equation}\label{Phi_Lag}
 \Phi(x, y, u, \alpha, \beta, \gamma):=\left[ 
\begin{array}{c}
     \nabla_{x, y}L_{\text{K}}(x, y, u, \alpha, \beta, \gamma)\\[0.75ex]
     \left(\phi_{\text{FB}} (\alpha_j, \; -G_j(x))\right)_{j=1, \ldots, p}\\[0.75ex]
     \left(\phi_{\text{M}}\left(u_i,\;-g_i(x,y),\; \nabla_y g_i(x,y)\gamma, \;\beta_i\right)\right)_{i=1, \ldots, q}
\end{array}
\right] =0.   
\end{equation}
A careful analysis of a nonsmooth Newton method to solve this system is conducted in \cite{harder2021reformulation}. Different approaches to construct numerical methods to directly compute different types of statioarity concepts for problem \eqref{eq:KKT-Optimistic-Bilevel} can be found in \cite{guo2015solving,wu2015inexact}. 

A typical challenges for a method to directly compute stationary points for problem \eqref{eq:KKT-Optimistic-Bilevel} via  a nonsmooth system of equations reside in the selection of transformation functions (such as $\phi_{\text{FB}}$ and $\phi_{\text{M}}$ above) and corresponding adequate generalized differentiation objects for Newton or Levenberg–Marquardt step; see Section \ref{sec:Some final thoughts} for a discussion on this topic and some potential ideas on how the current bilevel learning machinery can be used in this context to scale the corresponding techniques up. 

\subsection{The Big--M strategy} It is one of the most common approaches to solve problem \eqref{eq:KKT-Optimistic-Bilevel} in practice. It consists of replacing the product term $u^\top g(x,y)$ in the complementarity constraints of problem \eqref{eq:KKT-Optimistic-Bilevel} by the two conditions
\begin{equation}\label{Big-M-strategy}
    u_j \leq v_jM_D, \;\; -g_j(x, y)\leq (1-v_j)M_P, \;\;, v_j\in \{0, \, 1\},\; j=1, \ldots, q,
\end{equation}
where $M_D>0$ and $M_P>0$ are constants assumed to be \textit{large enough}; hence, they are called \textit{Big-M}s. This approach is typically used for fully linear bilevel optimization problems; i.e., the problem \eqref{eq:Bilevel-Optimization-Problem} where all the functions involved are linear in $(x, y)$. Observe that for such a problem, replacing $u^\top g(x,y)$ in \eqref{eq:KKT-Optimistic-Bilevel} with the system \eqref{Big-M-strategy} will lead to a linear problem, which is parameterized by the big-Ms. In this case, the resulting problem would be a linear program, which can then be embedded in an algorithmic process requiring standard off-the-shelf tools for linear programs. 
The first main challenge with this approach is that the big-Ms need to be chosen such that an optimal solution to  \eqref{eq:KKT-Optimistic-Bilevel} is not cut off. 

In \cite{kleinert2020there}, it is shown that identifying a suitable value for the big-M is an NP-hard problem. Moreover, even when a suitable value for the big--M could be found, solving the resulting optimization problem with constraint of the form is mixed-integer problem, which would not be scalable in the bilevel learning context. 
Note that an alternative to address the challenge with identifying a suitable big-M is to use an SOS1 scheme to construct a different constraint system, leading to a mixed-integer optimization that is much easier to solve; see, \cite{kleinert2023there} for a detailed exposition of the SOS1 scheme and its advantages of the big--M strategy. 

\subsection{KKT reformulation in bilevel learning}\label{KKT reformulation in bilevel learning}
The reformulation \eqref{eq:KKT-Optimistic-Bilevel} was the primary approach in the series of papers \cite{bennett2006model,bennett2006interplay,Kunapuli2007BilevelMS,Bennett2008ABO,bennett2008bilevel,Kunapuli2008ClassificationMS,moore2009nonsmooth,Bennett2010BilevelPA} by  Kristin Bennett and her co-authors,  focused on special version of the hyperparameter optimization problem \eqref{eq:Optimistic-Concrete}--\eqref{eq:lower-level-concrete}, especially support vector problems with linear kernel. These papers played a key role in the promotion of applications of bilevel optimization in machine learning, as already discussed earlier in Section \ref{sec:A short history}. The approach used in all these papers is based on first transforming the corresponding problems into the form \eqref{eq:KKT-Optimistic-Bilevel} and then applying off--the-shelf solvers on them. More recently, we have had many other papers studying hyperparameters optimization problems for linear or nonlinear support vector machines \cite{li2022bilevel,qian2023global,li2025constraint,coniglio2023bilevel,ward2025mathematical}, where the transformation \eqref{eq:KKT-Optimistic-Bilevel} is also the base for the numerical methods. In many of these papers (see, e.g., \cite{coniglio2023bilevel,Kunapuli2007BilevelMS,Kunapuli2008ClassificationMS, moore2009nonsmooth,qian2023global,ward2025mathematical,li2022bilevel}), it is shown that the BO approach can lead to algorithms that are more efficient than classical methods to conduct this process, such as grid search and Bayesian optimization, if these techniques are programmed under the same framework.  Of course, no complexity analysis has been conducted on such methods. However, the message that we can draw from the mentioned references is that if the lower-level problem is constrained the KKT reformulation discussed here has the potential to lead to efficient methods for BL. 

To close this section, observe that if the lower-level problem in \eqref{eq:Bilevel-Optimization-Problem} is unconstrained, the problem \eqref{eq:KKT-Optimistic-Bilevel} reduces to a problem with the constraint  $\nabla_y f(x, y)=0$, and hence no complementarity constraints appear in this case. In the papers \cite{alcantara2021unified,okuno2021}, this approach is used to deal with hyperparameter optimization problems where the lower-level problem is unconstrained (or with the lower-level constraints embedded to the lower-level objective function). However, as the lower-level objective function is nonsmooth there, the approach is extended by means a smoothing technique that replaces this objective function by an approximation function that enable the recovery of the some useful information (e.g., element from the the generalized Jacobian) from the original nonsmooth function when the involved parameter is driven to zero. In the papers \cite{alcantara2021unified,okuno2021}, it is also shown that the smoothing algorithm developed there is far much faster than grid search and Bayesian optimization methods for hyperparameter optimization.

\section{Lower-level value function reformulation-based methods}\label{sec:Pure first and second order methods}

One of the common points between the implicit function  model \eqref{eq:Implicit-Function-Model}, used in classical bilevel learning algorithmic framework  presented in Section \ref{sec:Main algorithmic techniques}, and the KKT reformulation \eqref{eq:KKT-Optimistic-Bilevel} covered in the previous section is that they both require second (resp. third) order derivative information on the lower-level problem for the development of  \textit{first} (resp. \textit{second}) order methods. 
    This higher order derivative requirement can be avoided if the model
\begin{equation}\label{eq:LLVF-Reformulation-Pso}\tag{LLVF}
      \underset{x, y}{\min}~F(x, y)\;\;  \mbox{s.t.}\;\, x\in X, \;\,  y\in Y(x),\;\, f(x,y)-\varphi(x)\leq 0,
\end{equation}
known as the lower-level value function (LLVF) reformulation, 
is considered as single-level transformation in the context of the standard optimistic problem  \eqref{eq:Standard-Optimistic-Bilevel}. Note that here, $\varphi$ represents the lower-level (optimal) value function 
\begin{equation}\label{eq:LLVF-varphi}
    \varphi(x):=\underset{y\in Y(x)}\min~f(x,y).
\end{equation}
It is important to observe that unlike for the implicit function model and the KKT reformulation, where lower-level convexity and constraint qualifications are needed in the transformation process of problem \eqref{eq:Standard-Optimistic-Bilevel}, no assumption is required to write problem \eqref{eq:LLVF-Reformulation-Pso}. This is due to the fact that by definition, it holds that 
\[
S(x)=\left\{\left. y\in Y(x)\right|~f(x,y)\leq \varphi(x)\right\}.
\]
Therefore, problems \eqref{eq:Standard-Optimistic-Bilevel} and \eqref{eq:LLVF-Reformulation-Pso} are globally and locally equivalent for free.

Before embarking on further in-depth analysis of problem \eqref{eq:LLVF-Reformulation-Pso}, let us place the approach in the historical context. Over 30 years ago, the problem \eqref{eq:LLVF-Reformulation-Pso} was considered in the literature from at least three different perspectives. First, Outrata \cite{outrata1988note,outrata1990numerical} explored the approach as a framework to numerically solve problem \eqref{eq:Standard-Optimistic-Bilevel}. Around the same time, Loridan and Morgan (see, e.g., \cite{morgan1989constrained,loridan1989new}) used \eqref{eq:LLVF-Reformulation-Pso} to establish regularity and stability results, where most often, the aim is to relax the value function constraint
$
f(x, y) - \varphi(x)\leq \epsilon
$
(with regularization parameter $\epsilon >0$) 
and establish some sequential convergence properties in relation to optimal solutions and value functions associated to the original optimistic and pessimistic models \eqref{eq:Optimistic-Model} and \eqref{eq:Pessimistic-Model}, respectively. Work on sequential stability analysis of this type has continued to be very  active (see, e.g., \cite{caruso2020regularization,caruso2025lower}). Subsequently, around the early 1990s, Ye and Zhu \cite{ye1995optimality} introduced the study of necessary optimality conditions for problem \eqref{eq:Standard-Optimistic-Bilevel} based on the LLVF reformulation. 
This is the stream of work that has mostly come to prominence in the last 30 years in relation to problem \eqref{eq:LLVF-Reformulation-Pso}. 

\subsection{Differentiability of the lower-level optimal value function}
We start here by considering the lower-level problem \eqref{eq:Lower-level-problem} while letting $U$ be an open neighborhood of a point $\bar x \in X$. Then the following implication holds:
\begin{equation}\label{eq:Implication_y()_phi}
    y(\cdot)\in \mathcal{C}^1(U) \quad \Longrightarrow \quad   \varphi \in \mathcal{C}^1(U).
\end{equation}
In fact, if $y(\cdot)\in \mathcal{C}^1(U)$, then for all $x\in U$, $\varphi(x)=f(x, y(x))$ and therefore, similarly to the  gradient formula for $\mathcal{F}$ in Section \ref{sec:Main algorithmic techniques}, we have 
\begin{equation}\label{eq:Nabla_Varphi_1}
    \nabla \varphi(x) = \nabla_x f(x, y(x)) + \nabla y(x)^\top \nabla_y f(x, y(x)). 
\end{equation}
Obviously, if the lower-level problem is unconstrained, it follows from the first order optimality conditions $\nabla_y f(x,y(x)) = 0$, that for all $x\in U$,
\begin{equation}\label{eq:Nabla_Varphi_0}
    \nabla \varphi(x) = \nabla_x f(x, y(x)).
\end{equation}
Otherwise, if the lower-lower feasible set-valued mapping $Y$ is defined as in \eqref{eq:Y(x)_g(x,y)}, i.e., $Y(x):=\left\{y\in \mathbb{R}^m\left|\;\, g(x, y)\leq 0\right.\right\}$, then it results from the Lagrange multiplier rule (see the left-to-right implication in \eqref{eq:LL_KKT_System}) that for any $x\in U$, under a lower-level constraint qualification such as the MFCQ, for example, at the point $(x, y(x))$, then we have
\begin{equation}\label{eq:Nabla_y_f()_1} 
  \nabla_y f(x, y(x)) + \nabla_y g(x, y(x))^\top u(x)=0.  
\end{equation}

On the other hand, with $I\equiv I(\bar x, y(\bar x))$, it follows from (A1)--(A4) that for some open neighborhood $U_0 \subset U$ of $\bar x$, we have
\[
g_i(x, y(x)) =0 \quad \mbox{ for } \quad i\in I, \;\, x\in U_0.
\]
Hence, by applying the chain rule once again,
\begin{equation}\label{eq:Nabla_x_g()_2}
    0=\nabla_x g_i(x, y(x))+\nabla y(x)^\top \nabla_y g_i(x, y(x)) \quad \mbox{ for } \quad i\in I, \;\, x\in U_0.
\end{equation}
It follows from a combination of  \eqref{eq:Nabla_Varphi_1}, \eqref{eq:Nabla_y_f()_1}, and \eqref{eq:Nabla_x_g()_2} that 
\begin{equation}\label{eq:Nabla_Varphi_3}
    \begin{array}{rll}
       \nabla \varphi(x) & = & \nabla_x f(x, y(x)) + \nabla y(x)^\top \left[-\sum_{i\in I} u_i(x) \nabla_y g_i(x, y(x)) \right],\\[2ex]
       & = & \nabla_x f(x, y(x))  -\underset{i\in I}{\sum} u_i(x)  \left[\nabla y(x)^\top\nabla_y g_i(x, y(x)) \right],\\[2ex]
       & = & \nabla_x f(x, y(x))  + \underset{i\in I}{\sum} u_i(x)\nabla_x g_i(x, y(x)).
    \end{array}
\end{equation}

Moreover, it is clear from \eqref{eq:Nabla_Varphi_0} and the last line of equation \eqref{eq:Nabla_Varphi_3} that when the value function $\varphi$ is smooth, unlike in the context of $y(\cdot)$, the expression of its gradient needs only first order information for the functions involved in the lower-level problem. This remains true even if $\varphi$ is not smooth. Before we discuss this aspect, note that the converse of implication \eqref{eq:Implication_y()_phi} is not true. In particular, for the example of parametric problem in \eqref{example-illus-1}, the optimal solution function $y(\cdot)$ is nonsmooth at $0$ and $1$ (see Figure \ref{fig:y() example}), while the value function $\varphi$ is continuously differentiable at these same points (see Figure \ref{fig:varphi_example}). 
\begin{figure}[H]
    \centering
    \includegraphics[scale = 0.45]{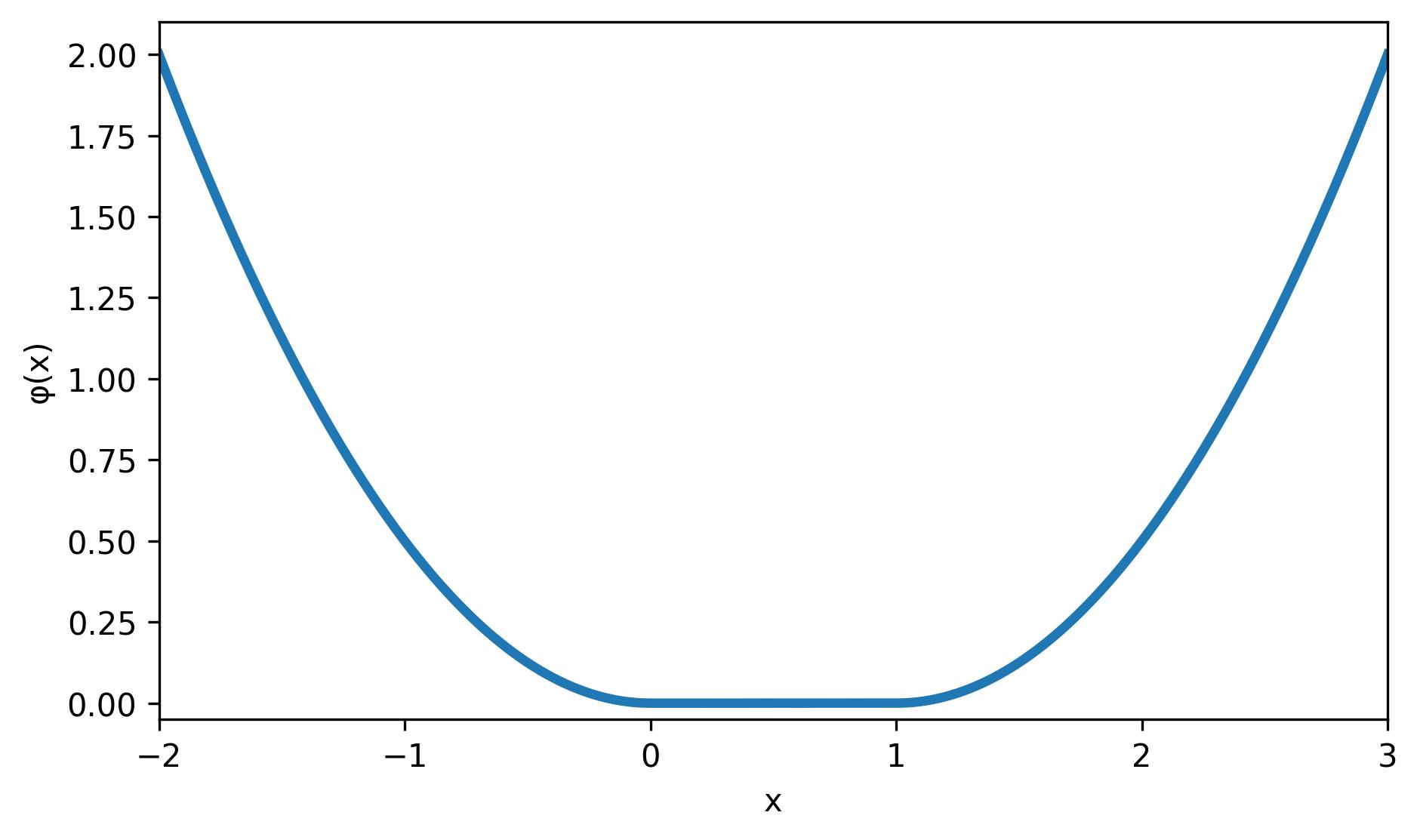}
    \caption{Graph of the lower-level value function $\varphi$ associated to problem \eqref{example-illus-1}.}
    \label{fig:varphi_example}
\end{figure}


Now, observe that even when $\varphi$ is nonsmooth, a subgradient of the function can still typically have the form in the last line of equation \eqref{eq:Nabla_Varphi_3}. To see this, assume that the lower-level feasible set-valued mapping $Y$ is given in \eqref{eq:Y(x)_g(x,y)} and let the functions $f$ and $g_i$ for $i=1, \ldots, q$ are smooth and fully convex; i.e., convex in the joint variable $(x, y)$. Furthermore, let $S(x) \neq \emptyset$ for all $x\in U \supset X$. Then $\varphi$ is locally Lipschitz continuous near any point $x\in X$; furthermore, for any $y\in S(x)$,  there exists $u\in \mathbb{R}^q$ such that 
\begin{equation}\label{eq:SubDiff}
 u\in \Lambda(x,y) \quad \mbox{ and }\quad   \nabla_x f(x, y) + \nabla_x g(x, y)^\top u \; \in \partial \varphi(x), 
\end{equation}
where $\partial \varphi(x)$ represents the subdifferential of $\varphi$ at $x$, in the sense of convex analysis. 
More precisely, if a constraint qualification (e.g., the MFCQ) holds at $(x, y)\in \mbox{gph}S$ for the lower-level constraint, then it holds that 
\begin{equation}\label{eq:Sub-varphi-conv}
    \partial \varphi(x) = \left\{\left.\nabla_x \ell(x, y, u)\right|~ u\in \Lambda(x, y)\right\};
\end{equation}
see, e.g., \cite{ye2023difference} for details on how to generate this formula. Recall that in the formula \eqref{eq:Sub-varphi-conv}, $\ell$ denotes the lower-level Lagrangian function \eqref{eq:Lower-level-Lagrangian}. 
In the case where the full convexity assumption on the lower-level problem is not satisfied, according to Gauvin and Dubeau \cite{gauvin2009differential}, if we assume that the set-valued mapping $Y$ \eqref{eq:Y(x)_g(x,y)} is nonempty and uniformly compact near $x$ and the MFCQ holds at $y$ in $Y(x)$ for fixed $x$ (for all $y\in S(x)$), then $\varphi$ is Lipschitz continuous near $x$ and its Clarke subdifferential can be estimated as
\begin{equation}\label{eq:Sub-Conv}   
  \partial \varphi(x) \subseteq \text{co}\,\left\{\underset{y\in S(x)}{\bigcup}\, \underset{u\in \Lambda(x, y)}{\bigcup} \left\{\nabla_x \ell(x, y, u)\right\} \right\},  
\end{equation}
where ``co'' stands for the convex hull of the corresponding set. This formula can be simplified in various ways, depending on the adjustments made on the assumptions; see, \cite{zemkoho2022estimates} and references therein for an overview of possible simplification scenarios. For instance, it naturally results from \eqref{eq:Sub-Conv} that $\varphi$ is strictly differentiable with 
\begin{equation}\label{Derivative-phi} 
    \nabla \varphi(x) = \nabla_x \ell(x, y, u) \;\mbox{ provided that }\; \{y\} = S(x) \mbox{ and } \{u\} = \Lambda (x, y).
\end{equation}

Thanks to the formulas of the subdifferential of $\varphi$ that have emerged from this discussion,  \textit{pure} first and second order methods can be constructed to solve problem \eqref{eq:LLVF-Reformulation-Pso}, as it will be clear below when we discuss solution algorithms. 
However, this does not necessarily make the problem easy to solve. For instance, similarly to problem \eqref{eq:KKT-Optimistic-Bilevel}, most classical CQs fail for problem \eqref{eq:LLVF-Reformulation-Pso}; see, e.g., \cite{dempe2011generalized,ye1995optimality}. Additionally, analogously to \eqref{eq:Implicit-Function-Model}, problem \eqref{eq:LLVF-Reformulation-Pso} is only implicitly defined and nonsmooth in general. 

\subsection{Optimality conditions}
Many papers have been written on  optimality conditions for problem \eqref{eq:LLVF-Reformulation-Pso}, mainly focusing on the use of variational analysis tool, considering the possible nonsmoothness of $\varphi$; see, e.g., \cite{ye1995optimality,dempe2007new,dempe2011generalized} and references therein. If $(x, y)$ is a local optimal solution of problem \eqref{eq:LLVF-Reformulation-Pso}, then, the possibly simplest class of necessary optimality conditions for problem \eqref{eq:LLVF-Reformulation-Pso} at this point is 
\begin{equation}\label{eq:StatsCond}\tag{StatsCond}
\left.\begin{array}{r}
\nabla_x F({x}, {y})+\nabla_x g({x}, {y})^\top (u - \lambda w)
+ \nabla G({x})^\top v =0 \\
 \nabla_y F({x}, {y}) +\nabla_y g({x}, {y})^\top (u-\lambda w) = 0 \\
\nabla_y f({x}, y) + \nabla_y g({x}, y)^\top w = 0\\
u\geq 0, \;\; g({x}, {y})\leq 0, \;\; u^\top g({x}, {y})=0\\
v\geq 0, \;\; G({x})\leq 0, \;\; v^\top G({x})=0\\
w\geq 0, \;\; g({x}, y)\leq 0, \;\; w^\top g({x}, y)=0\\[1ex]
\lambda\geq 0
\end{array}\right\}
\end{equation}
with $u$, $v$, and $\lambda$ representing upper-level Lagrange multipliers, respectively associated to the constraints $y\in Y(x)$, $x\in X$, and $f(x,y)-\varphi(x)\leq 0$, while using the description of the upper-level feasible set in \eqref{eq:X-G(x)}. Note that here, the vector $w$ corresponds to the lower-level Lagrange multiplier associated to the calculation of a subgradient of $\varphi$ at $x$ in the spirit of the formulas \eqref{eq:Sub-varphi-conv}, \eqref{eq:Sub-Conv}, and \eqref{Derivative-phi}.

An interesting feature of the stationarity system \eqref{eq:StatsCond} is that it does not explicitly depend on the optimal value function $\varphi$. However, feasibility is satisfied if the lower-level problem is convex, given that this ensures that $y\in S(x)$ if and only if  there exist $w\in \mathbb{R}^q$ such that the conditions in lines three and six are satisfied; cf. \eqref{eq:LL_KKT_System}.

The optimality conditions \eqref{eq:StatsCond} can hold at a point $(x, y)$ under the following assumptions (see \cite{ye1995optimality,dempe2007new,dempe2011generalized}, for example, for relevant theory):
\begin{itemize}
\item [(i)] Calmness of the set-valued mapping $\Phi :\mathbb{R} \rightrightarrows \mathbb{R}^m$ defined below at $(0, x, y)$:
\[
\Phi(\theta):=\left\{(x, y)\in \mathbb{R}^n\times \mathbb{R}^m\left|\; y\in Y(x), \;\, f(x,y)-\varphi(x) \leq \theta\right.\right\};
\]
    \item[(ii)] Fulfillment of the MFCQ at $x$ for the upper-level constraint described in \eqref{eq:X-G(x)};
    \item[(iii)] Fulfillment of the MFCQ at $y$ in $Y(x)$;
    \item[(iv)] Inner semicontinuous of $S$ at $(x, y)$ or full convexity of the lower-level problem.
\end{itemize}
Note that the lower-level problem is fully convex if the functions $f$ and $g_i$, for $i=1,\ldots, p$, are convex in $(x, y)$. The concept of inner semicontinuity of a set-valued mapping (see, e.g., \cite{mordukhovich2018variational}) is weaker than, but closely related to, the notion of lower semicontinuity, while, similarly, the calmness of a set-valued mapping is weaker than the Aubin property (extension of Lipschitz continuity to set-valued mappings). For the precise mathematical definition of  calmness and its characterizations and related applications, interested readers are referred to \cite{henrion2002calmness,henrion2005calmness}. For further discussion on all these assumptions in the context of bilevel optimization and the derivation of necessary optimality conditions for problem \eqref{eq:LLVF-Reformulation-Pso} interested readers can consult the references \cite{dempe2007new,dinh2010subdifferentials,mordukhovich2012variational,le2025regular,ye1995optimality,ye1997note,ye2004nondifferentiable,dempe2014necessary,dempe2011generalized,mehlitz2021note}.

For the avoidance of any doubt, it is worth recalling that necessary optimality conditions such as those in \eqref{eq:StatsCond} can be used to achieve at least three goals in relation to \eqref{eq:Standard-Optimistic-Bilevel}: (a) They can be used as stopping criteria for numerical algorithms to solve problem \eqref{eq:LLVF-Reformulation-Pso}. (b) They could be combined with suitable second order sufficient conditions to establish that a point is locally optimal  for  \eqref{eq:LLVF-Reformulation-Pso}; suitable second order sufficient conditions ensuring that necessary conditions of the type in \eqref{eq:StatsCond} can lead to locally optimal points for \eqref{eq:LLVF-Reformulation-Pso} are developed in \cite{fischer2022semismooth,mehlitz2021sufficient}. (c) Conditions such as \eqref{eq:StatsCond} can be used to build second order methods. More details on works related to points (a) and (c) in the development of numerical methods for \eqref{eq:LLVF-Reformulation-Pso} will be discussed later in this section. 


\subsection{Pure first order--type methods}
In \cite{outrata1988note}, a usual augmented Lagrangian function is considered and minimized with a bundle method; in this paper, we have unperturbed lower-level constraints, but no upper-level constraint. In \cite{outrata1990numerical}, the generalized equation, consisting of replacing $y\in S(x)$ by 
\[
    0\in \nabla_y f(x, y) +N_{Y(x)}(y),
\] 
where $N_{Y(x)}(y)$ denotes a normal cone to $Y(x)$ at the point $y$, is considered together with the implicit function and LLVF reformulations to address constrained lower-level optimistic bilevel programs, with special focus on the latter two models. The constraints of problem \eqref{eq:LLVF-Reformulation-Pso}  are fully penalized, under the assumption that the calmness condition, in the sense of Clarke \cite{clarke1990optimization}, is satisfied. 
For each of the transformations (i.e., the corresponding \eqref{eq:Implicit-Function-Model} and \eqref{eq:LLVF-Reformulation-Pso} problems), focus in \cite{outrata1990numerical} is on how to compute elements from the subdifferentials of the corresponding nonsmooth objective functions. This enables the use of the NDO (nondifferentiable optimization) solver based on the work in \cite{kiwiel1985methods} (see latest edition in \cite{kiwiel2006methods}) to compute solutions for the implicit function and LLVF models. 

 
Overall, in terms of developing first order methods, any algorithmic technique that requires only first order information to proceed can be explored in the context of problem \eqref{eq:LLVF-Reformulation-Pso}. Recently, a few papers on pure first order methods have appeared in the bilevel learning literature; see, e.g., \cite{liu2022bome,kwon2023fully,kwon2023penalty,gao2022value,ye2023difference}. The main approach in this context has consisted of solving the penalized problem 
\begin{equation}\label{eq:lambda-LLVF-Reformulation-Pso}
      \underset{x, y}{\min}~F(x, y) + \lambda \left(f(x,y)-\varphi(x)\right),
\end{equation}
provided there is no upper-- nor lower-level constraints. In this case, it obviously follows from \eqref{Derivative-phi} that $\nabla \varphi(x)=\nabla_x f(x, y)$ with $\{y\} = S(x)$, under suitable assumptions. In this context, a gradient descent scheme, similarly to Algorithm \ref{alg:1} for  problem \eqref{eq:Implicit-Function-Model}, can be developed for \eqref{eq:lambda-LLVF-Reformulation-Pso}.
One of the main challenges here 
identifying suitable values for the penalization parameter $\lambda$; see some relevant analysis in the paper \cite{tin2023levenberg}.

A notable representative of this line is the fully first-order stochastic approximation (F2SA) method of \cite{kwon2023fully}, which provides one of the cleanest complexity theories for the penalized formulation \eqref{eq:lambda-LLVF-Reformulation-Pso} in the unconstrained setting. The core difficulty is that the penalty term introduces a bias (since $\lambda<\infty$ does not enforce $y\in S(x)$ exactly), but taking $\lambda$ too large makes the penalized objective increasingly ill-conditioned, forcing smaller step sizes. F2SA resolves this trade-off by using a \emph{scheduled} penalty parameter $\{\lambda_k\}_{k\ge 0}$: starting from $\lambda_0>0$, it \emph{increases} $\lambda_k$ at a controlled polynomial rate (equivalently, it \emph{decreases} the effective penalty tolerance $1/\lambda_k$), while simultaneously shrinking the primal step sizes to maintain stability as the penalized landscape sharpens. This careful co-evolution of $(\lambda_k,\alpha_k,\beta_k)$ yields explicit non-asymptotic rates in both deterministic and stochastic regimes under the regularity assumptions used by implicit function methods, including iteration complexities to reach an $\epsilon$-stationary point scaling as $\tilde{\mathcal O}(\varepsilon^{-3/2})$ in deterministic settings and $\tilde{\mathcal O}(\varepsilon^{-5/2})$--$\tilde{\mathcal O}(\varepsilon^{-7/2})$ in stochastic settings depending on whether noise affects one or both levels. Closely related penalty-based value function methods, such as \cite{kwon2023penalty}, develop finite-time guarantees for first-order schemes by linking approximate stationarity of the penalized problem to approximate bilevel optimality through suitable choices of the penalty magnitude.

It is also instructive to contrast these lower-level value function penalization-based  rates with those obtained by the barrier-style value function method named BOME, published in \cite{liu2022bome}. While BOME likewise avoids implicit differentiation by exploiting the value-function constraint, its convergence guarantees are established for a \emph{different} notion of progress, namely a KKT-like residual tailored to the value-function constrained reformulation (rather than directly bounding the norm of the hypergradient.

When a general form of problem \eqref{eq:LLVF-Reformulation-Pso}, involving upper and/or lower-level constraints, is considered, penalizing only the value function constraint as done in \eqref{eq:lambda-LLVF-Reformulation-Pso}, corresponds to a partial penalization model, 
given that it remains constrained by both the upper- and lower-level constraints. 
This partial penalization approach was  introduced in \cite{ye1995optimality,ye1997note}, via a specialized version of the Clarke concept, label as \textit{partial calmness}, to derive necessary optimality conditions for problem \eqref{eq:LLVF-Reformulation-Pso}. The concept of partial calmness has since then become very prominent in the study of the LLVF reformulation and other optimization problem classes (see, e.g., \cite{ye1997exact,liu2008partial,zemkoho2021theoretical}). For instance, it can be used to replace the calmness concept in (i) above, as a qualification condition, to derive the necessary optimality conditions in \eqref{eq:StatsCond}. However, in this case, $\lambda$ will be positive and instead represent the penalty parameter, and not a Lagrange multiplier \cite{ye1995optimality,dempe2007new,dempe2011generalized,ye1998new}. 

Another approach that also builds on the partial penalization model in \eqref{eq:lambda-LLVF-Reformulation-Pso} is the difference-of-convex functions (DCA) method. To get a flavor of how this works, recall that if the lower-level objective and constraint functions in problem \eqref{eq:Standard-Optimistic-Bilevel} are fully convex, then the value function $\varphi$ \eqref{eq:LLVF-varphi} is convex. If we additionally assume that the upper-level objective function $F$ is convex in $(x, y)$, then the functions $F + \lambda f$ and $\varphi$ are both convex, and therefore, the objective function in \eqref{eq:lambda-LLVF-Reformulation-Pso} is the difference of two convex functions for a fixed $\lambda >0$.
Building on a well-established literature on DCA programming,  the paper \cite{gao2022value} studies this approach in the context of a bilevel hyperparameter selection problem; a particular feature of the DCA method, which differentiates it from  the gradient descent-type approaches in the previous papers (e.g., \cite{outrata1988note} or more generally in the bilevel learning literature), is that at each iteration, a linear approximation of the value function is built using an element from its subdifferential. Furthermore, in \cite{gao2022value}, the subproblem is made strongly convex by the addition of a proximal term. This approach is generalized in \cite{ye2023difference} to problems with DC upper-level objective function. The penalty methods in these two papers are same as in \cite{outrata1990numerical}, with only two differences: (i) the value function is approximated by a linear function based subgradient of $\varphi$ along the lines of the formulas \eqref{eq:Sub-varphi-conv}, \eqref{eq:Sub-Conv}, and \eqref{Derivative-phi}; and (ii) a proximal term is added to the subproblem for it to be strongly convex.

\subsection{Pure second order--type methods}
 Fixing $\lambda >0$ in the optimality conditions provided in \eqref{eq:StatsCond}, as it would be the case if the partial calmness concept is used as a qualification condition, the system can be written as 
\begin{equation}\label{Phi_lambda}
    \Phi_\lambda(x, y, u, v, w):= \left[\begin{array}{c}
         \nabla_{x,y}L_{\lambda}(x, y, u, v, w)\\[0.75ex]
         \nabla_y\ell(x, y, w)\\[0.75ex]
         \left(\phi_{\text{FB}} (u_i, \; -g_i(x, y)\right)_{j=1, \ldots, q}\\[0.75ex]
         \left(\phi_{\text{FB}} (v_j, \; -G_j(x))\right)_{j=1, \ldots, p}\\[0.75ex]
         \left(\phi_{\text{FB}} (w_i, \; -g_i(x, y)\right)_{j=1, \ldots, q} 
    \end{array}\right]=0,
\end{equation}
where the lower-level Lagrangian function $\ell$ and the Fischer-Burmeister function $\phi_{\text{FB}}$ are defined in \eqref{eq:Lower-level-Lagrangian} and \eqref{eq:phi_FB}, respectively. Note that the Lagrangian function $L_{\lambda}$ is defined for any $(x, y, u, v, w)$ with $(x, y)\in \mathbb{R}^n\times \mathbb{R}^m$ and $(u, w, v)\in \mathbb{R}^{2q}\times \mathbb{R}^p$ by
\[
L_{\lambda}(x, y, u, v, w):= F(x, y) + v^\top G(x) + (u-\lambda w)^\top g(x,y).
\]
Unlike the system \eqref{Phi_Lag} that involved second order information for the lower-level objective and constraint functions, \eqref{Phi_lambda} is a pure first order system. In fact, \eqref{Phi_lambda} is a $(n+m+p+2q) \times (n+2m+p+2q)$ semismooth system of equations. Based on these observations, Gauss-Newton and Levenberg--Marquardt methods are developed to solve the system in \cite{fliege2021gauss,tin2023levenberg,jolaoso2025fresh}. Since the system has $m$ more equations than variables, substituting $y$ in the second and fourth block of the right-hand-side of equation \eqref{Phi_lambda}  with a dummy variable $z$, we regularize \eqref{Phi_lambda} to a square semismooth system of equations, which can then enable the development of a semismooth Newton method; see \cite{zemkoho2021theoretical,fischer2022semismooth} for details in this direction, with a comparison of the numerical performance of \eqref{eq:KKT-Optimistic-Bilevel} and \eqref{eq:LLVF-Reformulation-Pso} conducted in \cite{zemkoho2021theoretical}. This comparison suggests that \eqref{eq:LLVF-Reformulation-Pso} generally leads to a better performance for the examples from the BOLIB library \cite{zhou2020bolib}, which are essenitally small toy problems. 

An important point that needs to be highlighted for the aforementioned works on solving the system \eqref{Phi_lambda} is that they are all \textit{pure} second order methods, as they do not require any third order derivative information. This would not be the case in the context of problems \eqref{eq:Implicit-Function-Model} and \eqref{eq:KKT-Optimistic-Bilevel} as such techniques will require third order derivative information for the functions involved in the corresponding lower-level problems. 

\subsection{Other types of LLVF--based methods}
Let us observe that the value function constraint, which represents the main component of the feasible set of problem \eqref{eq:LLVF-Reformulation-Pso}, can be equivalently written as 
\[
f(x, y)-f(x, z) \leq 0 \,\mbox{ for } \, z\in Y(x),
\]
which is a generalized semi-infinite constraint. 
Considering this representation, the series of papers 
\cite{mitsos2008global,djelassi2021recent,kleniati2014branch,wiesemann2013pessimistic} exploits 
semi-infinite type reformulations or related relaxation techniques  to build global optimization algorithms for problem \eqref{eq:Standard-Optimistic-Bilevel}, using branch and bound procedures as key ingredients to ensure numerical efficiency.  
The LLVF reformulation has also been used (see, e.g., \cite{tahernejad2020branch,fischetti2017new,fischetti2018use}) to develop cutting plane--type numerical algorithms for mixed-integer bilevel programs in the case where all involved functions are linear. 

The paper \cite{lin2014solving} proposes an approximation algorithm that relies on a concept of entropy integral function as a smoothing function for the optimal value function $\varphi$ when $Y(x)=Y$ (i.e., unperturbed). We also have an algorithm in \cite{dempe2014solution}, where the value function \eqref{eq:LLVF-varphi} for $f(x, y):=x^\top y$ and $Y(x):=\left\{y\in \mathbb{R}^m|\; Ay\leq b\right\}$ is iteratively approximated by a linear approximation. The advantage with the structure of the lower-level problem here is that as $Y$ is an unperturbed polyhedral set, an optimal solution for the lower-level problem can be found at one of its extreme points for all $x\in X$, under a mild assumption (e.g., if $Y$ is a bounded polyhedral). For a fully linear upper-level problem, it is shown that the algorithm converges to a global or local optimal solution, depending on the specific assumption scenario considered. The idea in \cite{dempe2014solution} is later extended to more general problem classes in \cite{dempe2016solution}.   However, similarly to most of the schemes just described above, the scalability of this class of method to relatively large problem classes is  uncertain. Hence, their applicability to bilevel learning problems might be very limited.

It might also be useful to mention that in \cite{lampariello2017bridge}, the  reformulation \eqref{eq:LLVF-Reformulation-Pso} is used as base for the construction of a Generalized Nash equilibrium problem (GNEP) that is closely related to problem \eqref{eq:Standard-Optimistic-Bilevel}. This GNEP model is then exploited in \cite{lampariello2020numerically} to build a numerical method to compute approximate stationarity points for problem \eqref{eq:LLVF-Reformulation-Pso}.

\section{Comparing the reformulations of the optimistic bilevel program}\label{sec:Comparing the reformulations}
The main focus of this paper is the standard optimistic bilevel optimization problem \eqref{eq:Standard-Optimistic-Bilevel}, expressed in the implicit function model \eqref{eq:Implicit-Function-Model} when condition \eqref{eq:S(x)=1} is satisfied, and otherwise in the KKT and LLVF reformulations \eqref{eq:KKT-Optimistic-Bilevel} and \eqref{eq:LLVF-Reformulation-Pso}, respectively, when \eqref{eq:S(x)>1} holds. 
It therefore seems natural to take a little moment in this section to briefly compare the single-level reformulations \eqref{eq:Implicit-Function-Model}, \eqref{eq:KKT-Optimistic-Bilevel}, and \eqref{eq:LLVF-Reformulation-Pso} of problem \eqref{eq:Standard-Optimistic-Bilevel}. 

Of course, these three problems are so significantly different from each other, with the challenging component of \eqref{eq:Implicit-Function-Model} appearing in its objective function, while problems \eqref{eq:KKT-Optimistic-Bilevel} and \eqref{eq:LLVF-Reformulation-Pso}  have very complex feasible sets. In terms of the smoothness of \eqref{eq:Implicit-Function-Model} and \eqref{eq:LLVF-Reformulation-Pso}, we have condition \eqref{eq:Implication_y()_phi}, which implies that under suitable conditions, \eqref{eq:LLVF-Reformulation-Pso} will be automatically smooth if \eqref{eq:Implicit-Function-Model} is. However, the converse of this implication is not true as demonstrated by the example in \eqref{example-illus-1}; cf. Figures \ref{fig:y() example} and \ref{fig:varphi_example}. For specific comparisons between \eqref{eq:KKT-Optimistic-Bilevel} and \eqref{eq:LLVF-Reformulation-Pso} from multiple perspectives, interested readers are referred to \cite{zemkoho2021theoretical}.

\newcommand{\pmark}{\raisebox{0.2ex}{\(\sim\)}} 

\renewcommand{\arraystretch}{1.2} 
\begin{table}[h!]\label{Table:Reformulations}
\centering
\caption{Some basic comparisons of the reformulations of the standard optimistic bilevel optimization problem. 
\pmark\ denotes partial guarantees (often for smoothed/MPEC surrogates, requiring bounded iterates).
$^{\ast}$ denotes rate results proved under lower-level strong convexity and lower-level Hessian smoothness (as in implicit-function analyses).
}
\begin{tabular}{ll||ccc}
\toprule
 & \textbf{Property} & \eqref{eq:Implicit-Function-Model} & \eqref{eq:KKT-Optimistic-Bilevel} & \eqref{eq:LLVF-Reformulation-Pso} \\
\midrule

\multirow{3}{*}{\textbf{SLR requirements}}
  & Smoothness      & \cmark & \cmark & \xmark \\
  & Convexity   & \cmark & \cmark & \xmark \\
  & Strong convexity   & \cmark & \xmark & \xmark \\
  & LLCQ        & \cmark & \cmark & \xmark \\
\midrule

\multirow{2}{*}{\textbf{ULCQ fulfilment}}
  & UMFCQ            & \cmark & \xmark & \xmark \\
  & Can UMFCQ be restored? & \cmark & \cmark & \xmark \\
\midrule

\multirow{2}{*}{\parbox[t]{3cm}{\textbf{Derivative}\\\textbf{requirement}\\\textbf{for methods}}}
  & 1D1OM & \xmark & \xmark & \cmark \\
  & 2D2OM & \xmark & \xmark & \cmark \\[2ex]
\midrule

\multirow{2}{*}{\textbf{Convergence}}
  & Deterministic rates        & \cmark & \pmark & \cmark$^{\ast}$ \\
  & Stochastic rates     & \cmark & \xmark & \cmark$^{\ast}$ \\
  
\bottomrule
\end{tabular}\label{Table:Reformulations}
\end{table}

Overall, Table \ref{Table:Reformulations} summarizes key comparisons between \eqref{eq:Implicit-Function-Model}, \eqref{eq:KKT-Optimistic-Bilevel}, and \eqref{eq:LLVF-Reformulation-Pso} from four perspectives (while assuming that the lower-level problem is constrained): 
(a) The requirements needed to formally write the corresponding single-level reformulation (SLR) of  \eqref{eq:Standard-Optimistic-Bilevel} with LLCQ standing  for the \textit{lower-level constraint qualification} in reference to whether one is needed to write the corresponding reformulation. 
(b) The behavior of each reformulation with regards to a suitable version of the MFCQ that we label here as upper-level constraint qualification (ULCQ); a key thing to note is that the MFCQ fails for both the KKT and LLVF reformulations, when their constraints are treated as usual equality and inequality constraints.  However, the MFCQ can be restored for  \eqref{eq:KKT-Optimistic-Bilevel} by suitably addressing the combinatorial structure in the complementarity conditions; see, e.g., \cite{dempe2012karush,scheel2000mathematical}. So far, no restoration approach for the MFCQ has been discovered for problem \eqref{eq:LLVF-Reformulation-Pso};  see, e.g., \cite{zemkoho2021theoretical} for a related discussion. 
(c) The third block of Table \ref{Table:Reformulations} corresponds to the derivative requirements for the corresponding reformulation in the sense that  the abbreviation 1D1OM is used to refer to \textit{whether only first order derivatives are enough to develop a first order method} for the corresponding reformulation, while 2D2OM refers to \textit{whether only second order derivatives are enough to develop  second order methods}. 

Finally, the fourth aspect (d) concerns the convergence guarantees currently available in the BL literature, distinguishing between deterministic and stochastic rates. From the standpoint of \emph{available complexity guarantees}, the picture is currently quite uneven across reformulations. 
For the LLVF  route, recent work has established explicit non-asymptotic rates for \emph{pure first-order} schemes, most prominently via barrier/penalty mechanisms; see BOME \cite{liu2022bome} and the penalty-based methods in \cite{kwon2023fully,kwon2023penalty}. These guarantees are typically developed for the \emph{unconstrained} bilevel setting (or after handling constraints separately) and require lower-level regularity such as strong convexity/PL-type conditions and smoothness to control the bias induced by finite penalty/barrier parameters and finite inner-loop accuracy. 
In contrast, for KKT-type reformulations, while classical theory clarifies how constraint qualifications may be recovered by treating complementarity carefully \cite{scheel2000mathematical,dempe2012karush}, finite-time rate results in modern machine learning settings are comparatively limited and often appear in forms that are \emph{partial} in the sense of Table~\ref{Table:Reformulations}: existing analyses typically provide decay bounds for a \emph{KKT residual} of a smoothed/regularized surrogate and rely on additional technical conditions such as bounded iterates/compactness and smoothing schedules (e.g., method-of-multipliers/augmented-Lagrangian developments); see, e.g., \cite{lu2022stochastic,liu2023averaged}.




\section{Conclusions and final remarks}\label{sec:Some final thoughts}
Based on the survey of the BL and BO literature conducted in this paper, the following concluding observations can be made:

(i) The implicit function approach, labelled as \eqref{eq:Implicit-Function-Model}, has been working very well in solving specific classes of BL problems, thanks to efficient approximations of the lower-level optimal solution function $y(\cdot)$ and its Jacobian $\nabla y(\cdot)$, when both functions are well-defined. However, this approach has many limitations, as not only, ensuring that the required basic assumption \eqref{eq:S(x)=1} is satisfied is very difficult, but making sure that the lower-level optimal solution function is smooth function is even harder; cf. discussion in Section \ref{sec:Challenges and limitations}. As also highlighted in the latter section, things get worse when the lower-level problem is contrained. 

(ii) For some special classes of BL problems, alternative methods to state of the art techniques, 
have been shown not only to better capture the corresponding task, but numerically,  their resulting BO formulation can be more efficiently solved, even when lower-level problem is constrained. This is the case, for example, for hyperparameter optimization is machine learning, as discussed in Subsection \ref{KKT reformulation in bilevel learning}. The BO formulation of the problem can be numerically solved more efficiently and accurately in comparison to the state of the art grid search and Bayesian approaches, which are standard in the main stream machine learning literature, and also the most widely used in practice. However, such BO--based tools have not yet made their way to main stream machine learning  infrastructures such as widely used open source libraries. 
A likely reason is not only limited awareness, but also the practical complexity of current BO solvers: despite their principled formulation and strong numerical results, many methods rely on several algorithmic \emph{meta-parameters}---for example, the number of inner iterations, truncation depth, damping/regularization, linear-solver tolerances, and step sizes for both upper- and lower-level updates. These choices can have a major effect on stability, memory footprint, and runtime, which makes robust off-the-shelf deployment in mainstream machine-learning libraries more difficult. The adaptive methods discussed in Subsection~\ref{sec:adaptive} might help mitigate these issues but are in early stage of development.

(iii) There is an exponential number of applications of bilevel optimization in machine learning. But for most of these problem--types,  state of the art BL methods cannot be applied. For example, we can mention problems with lower-level constraints,  where the lower-level optimal solution set-valued is not single-valued for some upper-level variables, as well problems with nonsmooth lower-level objective or constraint functions, as outlined in Section \ref{sec:A flavor of machine learning}. For many of these problems, like the pessimistic case that results from \eqref{eq:S(x)>1}, for example, not much progress has been made in solving them in the general BO literature. 

(iv) There is also a wide range of numerical techniques in the BO literature that remain  unexplored in the context of BL; this paper has provided a brief overview of these approaches, with the hope that in the near future, they will draw the attention that they deserve. As highlighted in Sections \ref{sec:Constrained optimization}--\ref{sec:Pure first and second order methods}, a key limitation of the classical BO algorithms is that many are hard to scale, especially when they involve lower-level constraints. The powerful derivative approximation schemes from state of the art numerical schemes for BL (cf. Section \ref{sec:Main algorithmic techniques}) are a potential way forward in the context of methods for problems \eqref{eq:KKT-Optimistic-Bilevel} and \eqref{eq:LLVF-Reformulation-Pso}. Work to deploy such ideas to solve the LLVF reformulation has started in the context of unconstrained lower-level problems as highlighted in Section \ref{sec:Pure first and second order methods}.

The second challenge that comes with inequality constraints is the combinatorial nature of the formula involved in BO numerical schemes (see, e.g., the combinatorial nature of the systems \eqref{eq:nablaY(x)} and \eqref{Phi_lambda}, as well as in the complementarity conditions in the feasible set of problem \eqref{eq:KKT-Optimistic-Bilevel}). With regards to this, there is a good chance that GPU--based methods to scale up optimization algorithms, as outlined in \cite{schubiger2020gpu,bishop2024relu}, for example, could be potential ways forward to enable classical BO methods to solve reasonable size BL problems. 
Additionally, the emergence of quantum computing and its potential to enable the scaling of enumeration--based algorithms (see, e.g., \cite{leenders2024integrating}) is also a path that, in combination with derivative approximation--based schemes, could provide paths to accelerate classical BO methods to solve realistic BL problems.  

\section*{Acknowledgments}\label{sec:Acknowledgments}
\href{https://sites.google.com/view/pradeepkumarsharma/home}{Pradeep Sharma} collected and organized some of the references that helped in the writing of Section \ref{sec:A flavor of machine learning}.
The fourth  author would like to thank Stefan G\"{u}ttel (The University of Manchester), Coralia Cartis (University of Oxford),  Harry Zheng and Panos Parpas (Imperial College), and Oliver Stein (Karlsruhe Institute of Technology) for their seminar invitations, where some of the material in this paper was presented, and the interesting discussions that  inspired the development of some of the elements in the text.

\bibliographystyle{siamplain}   

\end{document}